\documentclass[11pt]{amsart} 

\usepackage[all]{xy}
\usepackage{amsthm}
\usepackage{mathrsfs}
\usepackage{amssymb}
\usepackage{amsmath} 
\usepackage{array, url}
\usepackage{stmaryrd}
\usepackage{mathtools}
\usepackage{tensor}
\usepackage{stmaryrd}
\usepackage{pifont}
\usepackage{float}
\usepackage[labelformat=empty]{caption}
\usepackage[margin=1in]{geometry}

\let\oldmarginpar\marginpar
\renewcommand\marginpar[1]{\-\oldmarginpar[\raggedleft\footnotesize #1]%
{\raggedright\footnotesize #1}}

\newcounter{firstnumber}[section]
\newcounter{secondnumber}[firstnumber]
\newcounter{thirdnumber}[secondnumber]
\newcounter{fourthnumber}[thirdnumber]
\newcounter{fifthnumber}[fourthnumber]
\newcounter{currentdepth}

\renewcommand{\thefirstnumber}{\arabic{section}.\arabic{firstnumber}}
\renewcommand{\thesecondnumber}{\thefirstnumber.\arabic{secondnumber}}

\newcommand{\segment}[2]{%
    \setcounter{currentdepth}{1}%
    \def\thesubsection{\thefirstnumber}%
    \refstepcounter{firstnumber}\label{#1}
    \addtocounter{subsection}{-1}%
    \subsection{#2}}
\newcommand{\ssegment}[2]{%
    \setcounter{currentdepth}{2}%
    \def\thesubsection{\thesecondnumber}%
    \refstepcounter{secondnumber}\label{#1}
    \addtocounter{subsection}{-1}%
    \subsection{#2}}

\newenvironment{iquation}
 {%
  \ifnum\thecurrentdepth=1 \refstepcounter{secondnumber}
    \else
    \ifnum\thecurrentdepth=2 \refstepcounter{thirdnumber}
      \else
      \ifnum\thecurrentdepth=3 \refstepcounter{fourthnumber}
        \else
        \ifnum\thecurrentdepth=4 \refstepcounter{fifthnumber}
        \fi
      \fi
    \fi
  \fi
  \begin{equation}
 }
{
\end{equation}
}

\swapnumbers

\newcommand{\Rep}{\operatorname{\bf{Rep}}}

\newcommand{\from}{\leftarrow}
\newcommand{\xto}{\xrightarrow}
\newcommand{\xfrom}{\xleftarrow}
\newcommand{\surj}{\twoheadrightarrow}

\newcommand{\gr}{\mathrm{gr}\,}

\newcommand{\Spec}{\operatorname{Spec}}

\newcommand{\Pic}{\operatorname{Pic}}

\newcommand{\Li}{\operatorname{Li}}

\newcommand{\Hom}{\operatorname{Hom}}
\newcommand{\Ext}{\operatorname{Ext}}
\newcommand{\Isom}{\operatorname{Isom}}
\newcommand{\End}[1]{\mathrm{End} \, {#1}}

\newcommand{\Aut}{\mathrm{Aut}\,}

\newcommand{\Fil}{\operatorname{Fil}}
\newcommand{\ad}{\operatorname{ad}}
\newcommand{\Lie}{\operatorname{Lie}}

\newcommand{\Mod}{\operatorname{Mod}}

\newcommand{\m}[1]{\mathrm{#1}}
\newcommand{\fk}[1]{\mathfrak{#1}}

\newcommand{\bb}[1]{\mathbb{#1}}
\newcommand{\cl}[1]{\mathcal{#1}}

\newcommand{\la}{\lambda}
\newcommand{\ka}{\kappa}
\newcommand{\si}{\sigma}

\newcommand{\ga}{\gamma}
\newcommand{\al}{\alpha}

\newcommand{\om}{\omega}

\newcommand{\GL}{\operatorname{GL}}

\newcommand{\Gm}{{\mathbb{G}_m}}
\newcommand{\Qp}{{\QQ_p}}
\newcommand{\Zp}{{\ZZ_p}}

\newcommand{\V}{\bb{V}}
\newcommand{\ZZ}{\bb{Z}}
\newcommand{\GG}{\mathbb{G}}

\newcommand{\NN}{\bb{N}}
\newcommand{\QQ}{\bb{Q}}
\newcommand{\PP}{\bb{P}}
\newcommand{\UU}{\bb{U}}
\newcommand{\FF}{\bb{F}}
\renewcommand{\AA}{\bb{A}}

\newcommand{\Oo}{\mathcal{O}}

\newcommand{\Aa}{\mathcal{A}}
\newcommand{\Pp}{\mathcal{P}}

\newcommand{\inv}{^{-1}}
\newcommand{\til}{\widetilde}

\newcommand{\inj}{\hookrightarrow}

\newcommand{\iso}{\stackrel{\sim}{\to}}
\newcommand{\liso}{\;\stackrel{\sim}{\longrightarrow}\;}
\newcommand{\gen}[1]{\mathopen\langle#1\mathclose\rangle}

\newcommand{\ZSinv}{\ZZ[S \inv]}
\newcommand{\pZSinv}{(\ZSinv)}
\newcommand{\ZTinv}{\ZZ[T \inv]}
\newcommand{\pZTinv}{(\ZTinv)}

\newcommand{\mtEtQpf}{	{
				\tensor*[^{\m {mt}}]{\m{ \acute Et }}{^\Qp_{\m f}}}
												}
\newcommand{\Et}{\mtEtQpf}
\newcommand{\Etinv}{\Et\pZTinv}
\newcommand{\mtdRQp}{{\tensor*[^{\m{mt}}]{\dR}{}(\Qp)}}
\newcommand{\DR}{\mtdRQp}

\newcommand{\mtMotQ}	{
										{
										\tensor*[^{\m{mt}}]{\m{Mot}}{^\QQ}
										}
						}
\newcommand{\Mot}{\mtMotQ}
\newcommand{\Minv}{\mtMotQ\pZSinv}

\newcommand{\Zpn}{{\ZZ/p^n}}

\newcommand{\Zpstar}{\ZZ_p^*}

\newcommand{\mupn}{\mu_{p^n}}

\newcommand{\hh}{\fk{h}}

\newcommand{\tbp}{\vec{01}}

\newcommand{\dR}{{\rm {dR}}}
\newcommand{\et}{{\textrm {\'et}}}
\newcommand{\mot}{{\rm {mot}}}

\newcommand{\pioneinfty}{\pi_1^{[\infty]}}

\newcommand{\pione}[1]{\pi_1^{[#1]}}
\newcommand{\pioneinftyet}{\pi_{1,\textrm{\'et}}^{[\infty]}}

\newcommand{\pioneet}[1]{\pi_{1,\textrm{\'et}}^{[#1]}}
\newcommand{\pioneinftydr}{\pi_{1,\textrm{dR}}^{[\infty]}}

\newcommand{\pionedr}[1]{\pi_{1,\textrm{dR}}^{[#1]}}

\begin{document}

\SelectTips{cm}{11}

\title[the thrice punctured line in depth two]{Explicit Chabauty-Kim theory for the thrice punctured line in depth two}
\author{Ishai Dan-Cohen and Stefan Wewers}
\date{\today}
\maketitle

%%%%%%%%%%%%%%%%%%%%%%%%%%%%%%%

\begin{abstract}
Let $X= \mathbb{P}^1 \setminus \{0,1,\infty \}$, and let $S$ denote a finite set of prime numbers. In an article of 2005, Minhyong Kim gave a new proof of Siegel's theorem for $X$: the set $X(\mathbb{Z}[S^{-1}])$ of $S$-integral points of $X$ is finite. The proof relies on a `nonabelian' version of the classical Chabauty method. At its heart is a modular interpretation of unipotent $p$-adic Hodge theory, given by a tower of morphisms $h_n$ between certain $\mathbb{Q}_p$-varieties. We set out to obtain a better understanding of $h_2$. Its mysterious piece is a polynomial in $2|S|$ variables. Our main theorem states that this polynomial is quadratic, and gives a procedure for writing its coefficients in terms of $p$-adic logarithms and dilogarithms.

\vspace{3mm}


\end{abstract}

%%%%%%%%%%%%%%%%%%%%%%%%%%%%%%%
%%%%%%%%%%%%%%%%%%%%%%%%%%%%%%%
%%%%%%%%%%%%%%%%%%%%%%%%%%%%%%%
\section{Introduction}

\segment{120824-1}{}
Let $X = \PP^1\setminus\{0,1,\infty\}$ and let $S$ be a finite set of primes of $\ZZ$. Then the set $X(\ZSinv)$ is finite: indeed, this is a special case of Siegel's theorem. But the old proofs leave something to be desired. Siegel's original proof is not effective, as is the proof of Faltings' theorem. And while \textit{Baker's method} is effective in the sense that it produces a bound on the number of points, the resulting bound is usually very large. (Baker's method is also hard to grasp intuitively, at least for a geometer.) The problem, in any case, of actually computing the set of solutions remains open.

In 2005, Minhyong Kim published a new proof of the finiteness of $X(\ZSinv)$. His proof marked the birth of a new theory of integral points, and the beginning of a program which may lead to the development of an algorithm for finding integral points on hyperbolic curves over number fields. The program has recently been brought into focus by a conjecture formulated by its father, which may be regarded as a sort of compromise between certain aspects of the BSD conjecture and Grothendieck's section conjecture \cite{nabsd}.

Kim's construction is essentially a generalization of Chabauty's method to higher nonabelian levels. Indeed, it is Chabauty's method that is most often used \textit{in practice} to find rational points. Kim's higher nonabelian levels extend its range of application beyond its initially limited scope. As in Chabauty's method, we fix an auxiliary prime $p \notin S$, and we set out to study the $p$-adic structure of the set of integral points. There is then a tower
\[
\xymatrix{
\vdots \ar[d] 					& \vdots \ar[d]									\\
\m{Sel}_2 \ar[r]^{h_2} \ar[d] 	& \m{Alb}_2 \ar[d] 							\\
\m{Sel}_1 \ar[r]^{h_1} 			& \m{Alb}_1 								}
\]
of maps of affine, finite type $\Qp$-schemes, and for each $n$, a commuting square like so.
\[
\xymatrix{
X(\ZSinv) 		\ar[r] \ar[d]_\ka 		&
X(\Zp) 			\ar[d]^\al 			\\
\m{Sel}_n(\Qp) 		\ar[r]_{h_n}			&
\m{Alb}_n(\Qp)						}
\]
Here, $\m{Sel}_n$ is the \emph{Selmer variety} of Kim, while $\m{Alb}_n$ (in our special case) is simply the de Rham fundamental group of $X_\Qp$, whose definition will be recalled below. As for the morphisms, among several competing terminologies, we use \emph{(global) unipotent Kummer map} for $\ka$, \emph{unipotent Albanese map} for $\al$, and \emph{(global) unipotent $p$-adic Hodge morphism} for $h_n$. For precise definitions, about which the reader need not worry just now, we will ultimately have to refer to Kim's articles \cite{kimi}, \cite{kimii} (but see segment \ref{120823-1} below for a brief reminder).

The unipotent Albanese map is well understood. It is ``Coleman analytic" (regular functions on $\m{Alb}_n$ pull back to Coleman functions on $X(\Zp)$), and has dense image. A Coleman function is better than a general locally analytic funtion. Coleman functions are defined by iterated ``Coleman integrals". In the case of the thrice punctured line, they are generated by the $p$-adic polylogarithms of Furusho \cite{FurushoI}, \cite{FurushoII}. These, in turn, have local power series expansions analogous to those of complex polylogarithms.

By contrast, the unipotent $p$-adic Hodge morphism is not well understood. It follows nevertheless, from a dimension-count aided crucially by Soul\'e vanishing, that for $n$ sufficiently large compared to the size of $S$, $h_n$ cannot have dense image. This implies that there exists a regular function $f$ on $\m{Alb}_n$ vanishing on the image of $h_n$. Its pullback $\al^*f$ to $X(\Zp)$ will then be a nonconstant Coleman function which vanishes on $X(\ZSinv)$. Being locally analytic, such a function can have only finitely many zeros.

So goes Kim's proof of finiteness. Its potential to be effective depends on our ability to explicitly write down functions $f$ as above. There is reason to believe that the resulting bound on the number of points would be at least reasonable: as $n$ grows, the set $\al\inv(h_n(\m{Sel}_n(\Qp)))$ shrinks, and Kim's conjecture, in our special case, says that for $n$ sufficiently large,
\[
\al\inv(h_n(\m{Sel}_n(\Qp))) = X(\ZSinv) 
\,.
\]
The general statement, as well as a large collection of examples in which it can be verified, may be found in \cite{nabsd}. This motivates a careful study of the maps $h_n$.

The abelian case, given by $n=1$, reduces to the case $X = \Gm$. Of course, $\Gm(\ZSinv)$ will not be finite, but the same constructions apply, to give us varieties $\m{Sel}_n(\Gm)$, $\m{Alb}_n(\Gm)$, and maps $h_n(\Gm)$ between them. Then a straightforward application of well known results leads to the conclusions that for $n$ arbitrary, $\m{Sel}_n(\Gm) = \AA_\Qp^S$, $\m{Alb}_n(\Gm) = \AA_\Qp^1$, and $h_n(\Gm)$ is the linear map given by
\[
\tag{$*$}
h_n(\Gm)(x_1, \dots, x_s) = (\log q_1)x_1 + \cdots + (\log q_s)x_s
\,,
\]
where $S = \{q_1, \dots, q_s \}$, and $\log$ denotes the $p$-adic logarithm. We then have $\m{Sel}_1(X) = \m{Sel}_1(\Gm)^2$, $\m{Alb}_1(X) = \m{Alb}_1(\Gm)^2$, and $h_1(X) = h_1(\Gm)^2$.

We turn to the case $n=2$. It follows from Soul\'e vanishing that the map $\m{Sel}_2 \to \m{Sel}_1$ is an isomorphism. It is also well known that $\m{Alb}_2 = H$ is the Heisenberg group of $3\times 3$ upper triangular matrices with $1$'s on the diagonal, with the map $\m{Alb}_2 \to \m{Alb}_1$ being the projection onto the abelianization of $H$. So $h_2 = (h_{1,2}, h_{2,3}, h_{1,3})$ has three coordinates, of which the first two come from $h_1$.

It is $h_{1,3}$ that forms the first nontrivial piece of the tower, a polynomial in $2s$ variables, whose determination is the subject of this work. Our main results state that $h_{1,3}$ is bilinear (segment \ref{111212c}), and give a simple algorithm for writing its coefficients as polynomials in $(\log q_i)$'s and $(\Li q_i)$'s, where $\Li$ denotes the $p$-adic dilogarithm (segment \ref{111212a}). Both results are corollaries of a construction that revolves around the fundamental group of the category of mixed Tate motives, and capitalizes on generalities concerning mixed Tate categories. This culminates in the big diagram of segment \ref{111212c}. A simplified version, emphasizing its most salient features, appears in our outline below (segment \ref{120203-2}).

%%%%%%%%%%%%%%%%%%%%%%%%%%%
\segment{120823-1}{}
Central to Chabauty-Kim theory is the unipotent fundamental group in its various guises. As explained in \cite[\S15]{Deligne89}, the tangent vector
\[
\vec{01} = (\frac{\partial}{\partial t} \mbox{ at } 0)
\]
may serve as base-point. When working with the punctured line, the unipotent fundamental group $\pioneinfty(X,\vec{01})$ is first and foremost a pro-unipotent-group-object of the category of mixed Tate motives. We more often work with the quotients
\[
\pione{n} = \pioneinfty/Z^{n+1}
\]
along the descending central series, numbered so that $\pione{1} = (\pioneinfty)^\m{ab}$; their $p$-adic \'etale realizations are denoted $\pioneet{n}$. Let $T=S\cup \{p\}$, and let $G_T$ denote the Galois group of the maximal extension of $\QQ$ unramified outside of $T$. Then the moduli space of $G_T$-equivariant $\pioneet{n}$-torsors is in a natural sense a nonabelian cohomology variety, denoted $\m H^1(G_T, \pioneet{n})$. Torsors which become trivial upon base-change to $B_\m{cris}$ form a closed subscheme $\m H^1_\m{f}(G_T, \pioneet{n})$. This is Kim's Selmer variety, which was denoted by $\m{Sel}_n$ above. The unipotent Kummer map $\ka$ sends a point $x \in X(\ZSinv)$ to the torsor of homotopy classes of paths $\vec{01} \to x$ suitably interpreted. There is also a local variant $\m H^1_\m{f}(G_p, \pioneet{n})$ of the Selmer variety, in which $G_T$ is replaced by a decomposition group $G_p$ at $p$, as well as a corresponding local variant of the unipotent Kummer map, $\ka^\m{loc}: X(\Zp) \to \m H^1_\m{f}(G_p, \pioneet{n})$. The local and global Kummer maps are compatible with the natural localization map $l: \m H^1_\m{f}(G_T, \pioneet{n}) \to \m H^1_\m{f}(G_p, \pioneet{n})$.

The unipotent fundamental group and its quotients also have realizations in the category of admissible filtered $\varphi$-modules over $\Qp$, denoted $\pionedr{\infty}$, $\pionedr{n}$, which we refer to vaguely as the \emph{de Rham} fundamental group. It follows from the unipotent $p$-adic Hodge theory of Olsson \cite{OlssonTowards} that $\pionedr{\infty}$ doubles as the fundamental group of the Tannakian category of unipotent vector bundles with integrable connection on $X_\Qp$ at a fiber functor associated to $\vec{01}$. The unipotent Albanese map $\al: X(\Zp) \to \pionedr{n}(\Qp)$ was considered first by Furusho \cite{FurushoI}, \cite{FurushoII}.  The main ingredient in its construction is the unique Frobenius-invariant path $\vec{01} \to y$, for any $y \in X(\Zp)$, which was discovered by Besser \cite{Besser}, and independently by Vologodsky \cite{Vologodsky}. In the case of the punctured line, the de Rham torsor of homotopy classes of paths $\vec{01} \to y$ actually has an obvious trivialization in addition to the one given by the Frobenius-invariant path; dividing the two gives us a point $\al(y) \in \pionedr{n}$.

Given a $G_p$-equivariant $\pioneet{n}$-torsor $P$, we may apply the Dieudonn\'e functor $D$ of $p$-adic Hodge theory to its coordinate ring, to obtain a $\pionedr{n}$-torsor, $P^\dR$. (Here the word \textit{torsor} is to be understood in the sense of ``algebraic geometry in a Tannakian category" \cite[5.4]{Deligne89}.) As before, $P^\dR$ will have two competing trivializations. Let us recall these briefly, referring to \S1 of \cite{kimii} for details.

The coordinate ring $\cl{P}^\dR = \Oo (P^\dR)$  is endowed with a descending filtration
\[
\Pp^\dR = F^0 \Pp^\dR \supset F^1\Pp^\dR \supset F^2\Pp^\dR \cdots
\]
by ideals known as the \emph{Hodge filtration}. We set
\[
F^i P^\dR := \Spec \big( \Pp^\dR / F^{-i+1}\Pp^\dR \big)
\,.
\]
The $0^\m{th}$ piece $F^0 \pionedr{n}$ of the fundamental group is the subgroup which corresponds to the category $\m{Vect}^{[\infty]}(\bar X)$ of unipotent vector bundles on the associated complete curve, and $F^0 P^\dR$ is the $F^0 \pionedr{n}$-torsor of tannakian paths that factor though $\m{Vect}^{[\infty]}(\bar X)$ (c.f. Section 2.2 of Hadian \cite{Hadian}). Since in our case $\bar X = \PP^1$ has no nontrivial unipotent vector bundles, we have
\[
F^0 \pionedr{n} = \{1\}
\]
trivial, so that $F^0P^\dR$ contains one and only one $\Qp$-point $p^H$.

The coordinate ring $\Pp^\dR$ is also endowed with an automorphism $\phi$, \textit{Frobenius}. It follows from Theorem 3.1 of Besser \cite{Besser} that $\pionedr{n}$ possesses a unique $\Qp$-point $p^\m{cr}$ fixed by $\phi$. We define $h_n^\m{loc}(P) \in \pionedr{n}(\Qp)$ to be the unique point for which
\[
p^\m{cr} h_n^\m{loc}(P) = p^H
\,.
\]

The crystalline Dieudonn\'e functor is an equivalence of categories (from crystalline representations to weakly admissible filtered $\varphi$-modules) which, according to Olsson's theory, sends $\pioneet{n}$ to $\pionedr{n}$. Combined with Proposition 1 of Kim \cite{kimii}, this means that $h_n^\m{loc}$ is an isomorphism
\[
\m H^1_\m{f}(G_p, \pioneet{n}) \xto{\cong} \pionedr{n}
\,;
\]
that it interchanges $\ka^\m{loc}$ and $\al$ is another consequence of Olsson's theory. Although we've limited our discussion to $\Qp$-valued points, both the localization map $l$ and the local unipotent Hodge morphism $h_n^\m{loc}$ are in a natural sense maps of schemes. Their composite $h_n = h_n^\m{loc} \circ l$ is the global unipotent $p$-adic Hodge morphism.

The resulting diagram
\[
\xymatrix{
X(\ZSinv) 
\ar[rr] \ar[d]_\ka 
&& X(\Zp) 
\ar[dr]^\al \ar[d]_{\ka^\m{loc}}  \\
\m H^1_\m{f}(G_T, \pioneet{n}) 
\ar[rr]_l 
&& \m H^1_\m{f}(G_p, \pioneet{n}) 
\ar[r]_-{h_n^\m{loc}}
& \pionedr{n}
}
\]
serves as a powerful tool for cutting out the integral points inside $X(\Zp)$; we will refer to it as ``Kim's cutter".

%%%%%%%%%%%%%%%%%%%%%%%%%%%%%%
\segment{120203-2}{}%%%%%
Our main definitions and results, in rough outline, are as follows. Recall that a \emph{mixed Tate category} $T$ over a field $k$ consists of a Tannakian category, plus an object $k(1)$, such that the $k(i):= k(1)^{\otimes i}$ are pairwise nonisomorphic objects of rank one, such that every simple object is isomorphic to $k(i)$ for some $i \in \ZZ$, and finally, such that the Yonneda Ext groups $\Ext^1(k(0), k(n))$ vanish for $n\le 0$. Mixed Tate categories give rise to a notion of \emph{mixed Tate duality}, similar to, but simpler than, Tannaka duality. It establishes a correspondence between mixed Tate categories over $k$ and prounipotent groups $U$ over $k$ equipped with a grading of the associated Hopf algebra $A$ such that $A_0= k$ and $A_i=0$ for $i < 0$. In contrast to Tannaka duality, the correspondence does not depend on the choice of a fiber functor. Moreover, the coordinate ring $A$ admits an elementary construction (\S3) in terms of certain equivalence classes of \emph{framed} objects.

A mixed Tate category $T$ admits a \emph{Heisenberg object} $H_T$, which, in turn, has a nonabelian cohomology set $\m H^1(T, H_T)$ which parametrizes isomorphism classes of \emph{Heisenberg $T$-torsors} (\ref{111018a}). In this language, we construct a diagram like so (\S4).
\[
\xymatrix{
\Ext^1_T(k(0),k(2))			\ar[d]_-g 						&
\m H^1(T, k(2)) 				\ar@{=}[l] \ar[d] 				\\
A_2							\ar[d]_-\nu					&
\m H^1(T,H_T)					\ar[l]_-{\la_2}  \ar[d]			\\
A_1 \otimes A_1				\ar[d]_-D 						&
\m H^1(T, k(1) \oplus k(1)) 	\ar[l]_-{\la_1}  \ar[d]^-\Delta	\\
\Ext^2_T(k(0),k(2))											&
\m H^2(T, k(2))					\ar@{=}[l]						}
\]
In it, the left column forms an exact sequence of $k$-vector spaces, while the right column is an exact sequence of pointed sets.

We now begin to hone in on our concrete situation, by considering $\pione{2}$, the depth two quotient of the unipotent fundamental group of the thrice punctured line at the base point $\vec{01}$. We conclude (\S5) that $\pione{2}$ is simply the Heisenberg group object of the mixed Tate category of mixed Tate motives. This means in particular that if we denote its Lie algebra by $\fk p^{[2]}_1$, then the exact sequence
\[
1 \to \QQ(2) \to \fk p^{[2]}_1 \to \QQ(1)^2 \to 1
\]
of mixed Tate motives is canonically split.

We may simplify our situation by replacing $X$ by the multiplicative group (\S6). This causes $\pioneinfty$ to be replaced by $\QQ(1)$, so Kim's cutter (\ref{120823-1}) looks like so.
\[
\xymatrix{
\ZZ[S\inv]^*					\ar[r] \ar[d]_\ka					&
\ZZ_p^*						\ar[d]^\al						\\
\m H^1_\m f(G_T, \Qp(1))		\ar[r]_-h							&
\Qp(1)_\dR														}
\]
Here $\Qp(1)_\dR$ denotes the Tate object of the Tannakian category of admissible filtered $\varphi$-modules over $\Qp$. We find that $\m H^1_\m f(G_T, \Qp(1)) = \AA^S_\Qp$, and that $h$ is the linear map (\ref{120824-1}($*$)).

Section 7 is devoted to a discussion of the local situation. Its main purpose is to give an explicit formula for the unipotent Albanese map. The result here, which is a direct consequence of formulas obtained by Furusho \cite{FurushoI}, \cite{FurushoII}, is as follows. Since $\pionedr{2}$ admits a canonical pair of generators, we may identify it with the group $H_\Qp$ of upper triangular $3 \times 3$ matrices with $1$'s along the diagonal. Then $\al:X(\Zp) \to H_\Qp$ is given by
\[
z \mapsto
\left(
\begin{matrix}
1 & \log z & -\Li z \\
0 & 1 & \log (1-z) \\
0 & 0 & 1
\end{matrix}
\right) \;,
\]
where $\Li$ denotes the $p$-adic dilogarithm.

In the language of mixed Tate categories, Kim's cutter (\ref{120823-1}) for $n=2$ may be rewritten as follows (\S8).
\[
\xymatrix{
%1
X(\ZZ[S \inv])	\ar[r] \ar[d]_\ka				&
X(\Zp)			\ar[d]^\al					\\
\m{H}^1(H_\et)		\ar[r]_{h_2}					&
\m{H}^1(H_\dR)								}
\]
Here, $H_\et$ denotes the Heisenberg group object of the mixed Tate category $\Etinv$ over $\Qp$ of mixed Tate representations of $G_T$ crystalline at $p$, while $H_\dR$ denotes the Heisenberg group object of the mixed Tate category $\DR$ over $\Qp$ of mixed Tate filtered $\varphi$-modules, and $h_2$ is the operation on Heisenberg torsors induced by a morphism $\Etinv \to \DR$ of mixed Tate categories.

It is convenient (but not necessary) to add to the picture a mixed Tate category $\Minv$ over $\QQ$, equipped with realization functors to $\Etinv$ and to $\DR$, satisfying $\m H^1 (\QQ(2)) = \m H^2(\QQ(2)) = 0 $, satisfying $\m H^1 (\QQ(1)) = \ZSinv^* \otimes \QQ$, and such that $\ka$ factors through a morphism $\tau: X\pZSinv \to \m H^1(H_\mot)$ to the cohomology of the Heisenberg group object. Although the best known candidate for such a category, and, indeed, the one on which we focus in \S8 below, is the category of mixed Tate motives, there are other alternatives: for instance, a thing that might be described as the category of \emph{mixed Tate hyperplane arrangements}, constructed in a four way collaboration by Beilinson, Goncharov, Schechtman and Varschenko \cite{BGSV}.

The title of \S9 refers not to one, but to two splittings. We have, on the one hand, a projection $\Sigma$ of the injection $\Qp(2)_\dR \inj H_\dR$ constructed in section 5 on the level of mixed Tate motives. If we let $A^\dR$ denote the graded Hopf algebra  associated to $\DR$, then we have, on the other hand, a projection $\widetilde \Psi$ of the map $g^\dR: \Ext^1_\DR(\Qp(0), \Qp(2)) \to A^\dR_2$ constructed in section 4 in an arbitrary mixed Tate category. These are compatible in the sense that they fit into a commuting square like so.
\[
\xymatrix{
A_2^\m{dR}					\ar[r]^-{\widetilde \Psi}						&
\Ext^1_\DR(\Qp(0), \Qp(2))	\ar@{=}[r]									&
\Qp																			\\
\m H^1(H_\m{dR})				\ar[u]^{\la_2^\dR} \ar[rr]_-{\sim} ^-\Psi		&
																			&
H_\m{dR}						\ar[u]_\Sigma								}
\]
As a corollary, we obtain an isomorphism $A_2^\dR = \QQ_p^2$.

We now come to our main contribution: the construction of the following commutative diagram (\S10).
\begin{small}
\[
\xymatrix{
A_2^\mot					\ar[d]_-\wr \ar[r]^{\tensor*[^A]{\rho}{^\dR_2}}	&
A_2^\dR					\ar[r]^-{\widetilde \Psi}							&
\Qp(2)_\dR																	\\
(A_1^\mot)^{\otimes 2}													&
\m H^1(H_\et)				\ar[r] \ar[d]^-\wr								&
H_\dR						\ar[u]_\Sigma									\\
\m H^1(\QQ(1)^2_\mot)	\ar[r] \ar[u]										&	
\m H^1(\Qp(1)_\et^2)		}
\]
\end{small}
This gives us an alternative, linear route, from the Selmer variety to the third de Rham coordinate axis, thus circumventing the structureless $\m H^1(H_\et)$ and dodging the ``Heisenberg coboundary equation" (see below). We let $\tensor*[_\QQ]{h}{_{1,3}}$ denote the resulting map $\m H^1(\QQ(1)^2) \to \Qp(2)_\dR$. We obtain the following

\subsection*{Theorem}
The third coordinate $h_{1,3}$ of the global unipotent $p$-adic Hodge morphism for the thrice punctured line in depth two over $\ZZ[S \inv]$ is bilinear.

\subsection*{}
(\ref{111212c}). We also obtain two formulas: (\ref{111227c}) the \textit{geometric values} for $t \in X\pZSinv$ of $\tensor*[_\QQ]{h}{_{1,3}}$ are given by
\[
\tensor*[_\QQ]{h}{_{1,3}}(t, 1-t) = -\Li(t)
\,,
\]
and (\ref{111227d}, the ``twisted anti-symmetry relation") for $u, v \in \Gm\pZSinv$, we have
\[
_\QQ h_{1,3}(u, v) + \tensor*[_\QQ] {h}{_{1,3}}(v,u) = (\log u)(\log v)
\,.
\]

Let $E = \QQ \otimes \QQ^*$, which we regard as a $\QQ$-vector space, written multiplicatively. Section 11 begins with an algorithmic proof of the vanishing $K_2(\QQ) \otimes \QQ = 0$ based on Tate's computation of $K_2(\QQ)$. This amounts to an algorithm for decomposing an arbitrary generator $q \otimes q'$ of $E \otimes E$ as a linear combination like so.
\[
q \otimes q' = \prod_k \big((u_k \otimes v_k) \cdot (v_k \otimes u_k)\big)^{s_k} \cdot \prod_l \big(t_l\otimes (1-t_l)\big)^{d_l}
\]
Here, $s_k$, $d_l$, are rational coefficients, and for each $k$ and $l$, $u_k$, $v_k$, $t_l$ are vectors $\in \QQ^*$.

In terms of this decomposition, we obtain a simple formula for the coefficients of $_\QQ h_{1,3}$ (or, equivalently, of $h_{1,3} $). For any two primes $q, q'$ distinct from $p$, we let $\tensor*[_\QQ ]{ h}{_{1,3}^{q,q'}}$ denote the $p$-adic number given by
\[
\tensor*[_\QQ ] {h}{^{q,q'}_{1,3}} = \sum_k s_k (\log u_k)(\log v_k) -\sum_l d_l \Li(t_l) \,.
\]

\subsection*{Theorem}
For any pair $x$, $y$, of elements of $\QQ^S$, we have the formula
\[
_\QQ h_{1,3}(x, y)
=
\sum_{1 \le i,j \le s} \tensor*[_\QQ ] {h}{^{q_i, q_j}_{1,3}}x_i y_j
\,,
\]
where $q_1, \dots, q_s$ denote the primes of $S$ (segment \ref{111212a}). 

\subsection*{}
We also obtain a formula for the map $\tensor*[^A]{\rho}{^\dR_2}: A_2^\mot \to A_2^\dR$: through the isomorphisms $A_2^\mot = \QQ^S \otimes \QQ^S$ and $A_2^\dR = \Qp \times \Qp$, it is given by
\[
\tensor*[^A]{\rho}{^\dR_2} (x \otimes y) = \big( \sum_{i,j=1}^s (\log q_i)(\log q_j) x_iy_j , \sum_{i,j=1}^s \tensor*[_\QQ ] {h}{_{1,3}^{q_i, q_j}} x_i y_j  \big)
\,.
\]
Actually, for our proofs of these formulas to work, we have to ensure that the rational numbers $u_k$, $v_k$, $t_k$ occurring in the decomposition of each $q_i \otimes q_j$ are coprime to $p$. We deal with this problem by assuming that $p$ is greater than the primes $q_i \in S$.

One noteworthy feature of these formulas is their independence of $p$. This is in a way explained by the motivic origin of $h$. Apropos \textit{motivic}, we should mention that our motivic relative $\m H^1(H_\mot)$ of the Selmer variety (which plays a role, for instance, in showing that our main diagram commutes) originates from the thesis of Majid Hadian-Jazi \cite{Hadian}.

The number $2$ being rather small, concrete applications of our results are not surprisingly limited. Coleman observed in his 1982 article \cite{Coleman} that the function
\[
D(z) = 2 \Li (z) - \log(z)\log(1-z)
\]
on $X(\Zp)$ for any odd $p$, vanishes on $X(\ZZ[1/2])$. As one application, we show that Coleman's observation squares with Chabauty-Kim theory. The result is the following proposition  (\ref{Examples1}), which may also be proved by less highfalutin arguments.

\medskip
\noindent
\textit{
Let $R$ be either $\ZZ[\ell^{-1}]$, for a prime
$\ell$, or $R=\mathcal{O}_K$ for a real quadratic number field $K$. Suppose,
moreover, that $p=11$ splits in $K$. Then the set $X(R)$ is empty, except if
$R=\ZZ[1/2]$ or $R=\ZZ[(1+\sqrt{5})/2]$, in which case we have
  \[
       X(\ZZ[1/2]) = \{-1,2,1/2\}
  \]
  or
  \[
       X(\ZZ[(1+\sqrt{5})/2]) = \{\frac{1\pm\sqrt{5}}{2},
            \frac{-1\pm\sqrt{5}}{2},\frac{3\pm\sqrt{5}}{2}\}.
  \]
}

\medskip
\noindent
This example provides a modicum of evidence for Kim's conjecture; see \S7 of \cite{nabsd}.

As a second application we obtain the following result (\ref{Examples2}), which is not new.

\medskip
\noindent
\textit{
Consider the exponential diophantine equation
\[
a\cdot b^x - c \cdot d^y = 1
\]
with $a,b,c,d \in \QQ^*$ fixed. Then the set of solutions $(x,y) \in \ZZ^2$ is finite.
}

\medskip
\noindent
Of course, we are more interested in obtaining a bound on the number of points. For instance, for the particular equation
\[
7^x - 3\cdot 2^y = 1
\,,
\]
our methods show that there are at most five solutions; in fact, there are only two.

%qqq

This article used to have an appendix. In it, we recounted our attempt to proceed by disrobing Kim's original construction. Starting with an element of 
\[
\m H^1_\m f (G_T, \Qp(1)^2)
\,,
\]
we have to make its unique lifting to $\m H^1_\m f (G_T, H_\et)$ in some way explicit. This gives rise to what we call the \emph{Heisenberg coboundary equation}, a problem in nonabelian Galois cohomology. Our analysis gives evidence for the advantage of our indirect approach, while being of potential interest in its own right. The interested reader will find it online \cite{HCE}.

\segment{120205-1}{}%%%%%%%%%%
The question posed in this article may equally well be posed for an arbitrary hyperbolic curve over an arbitrary number field in depth $n$. So here we are restricting attention to a special case in three ways. In one way this restriction is spurious: $\QQ$ may be replaced throughout by a totally real number field, at least for the results concerning the motivic version $_\QQ h_{1,3}$ of the object of interest. There is one exception: the algorithmic vanishing of $K_2 \otimes \QQ$ would be far more complicated.

In the other two directions, things get significantly more complicated immediately. Nevertheless, we hope to pursue both directions, starting as soon as we've finished writing the acknowledgements.

\subsection*
%{120205-2}
{Acknowledgements}
During the writing of this article, the first author benefited from meetings with Aravind Asok in Los Angeles, and with Jordan Ellenberg and David Brown in Madison; he is grateful for their hospitality, as well as for helpful conversations about the present work. Jordan Ellenberg, in particular, suggested a formula which is similar to our Theorem 11.2. The first author is also grateful to Alexander Goncharov for a helpful exchange of emails.

The second author's interest in this topic was raised during an invitation to the Renyi Institut in Budapest and from countless stimulating discussions with his host, Tamas Szamuely, in Budapest and in many other places. 

Both of us are grateful to Minhyong Kim for helpful conversations and for his encouragement to work on the problem implied by the title of this paper. Finally, we wish to thank the referee for his or her careful reading and for many helpful comments.

%%%%%%%%%%%%%%%%%%%%%%%%%%%%%%%%%
\section{Heisenberg extensions and the cup product}
\label{1110181}%%%%%
%%%%%%%%%%%%%%%%%%%%%%%%%%%%%%%%%

\segment{110522b}{}{
We fix a group $G$, and a ring $A$ containing $1/2$, and we consider representations
$
\rho_1, \rho_2: G \rightrightarrows \Aut E_1, E_2
$
over $A$. We let
$$
\rho_3= \rho_1 \otimes \rho_2: G \to \Aut E_3=E_1 \otimes_A E_2 \,.
$$
Here each $E_i$ is an $A$-module. For our present purpose $A$ may harmlessly be assumed to be a field, and each $E_i$ to be free of rank one.
}
%\segment{110522d}{}
More generally, $G$ may be a topological group, and $A$ a topological ring, in which case all representations, cohomologies, cocycles, etc., are to be assumed continuous.

\ssegment{110522e}{Definition}
We define the \emph{Heisenberg Lie algebra extension}, $\hh$, as follows. As a representation, we set
$
\hh := E_1 \oplus E_2 \oplus E_3 \,.
$
Given $x,y \in \hh$, we set
$
[x,y]:= (0,0, x_1\otimes y_2 - y_1\otimes x_2) \;.
$
This gives $\hh$ the structure of a Lie algebra, nilpotent of depth two. Its Baker-Campbell-Hausdorff product is thus given by
\begin{align*}
x\star y 	&= x+y+\frac{1}{2}[x,y] \\
			&= (x_1+y_1, x_2+y_2, x_3+y_3+\frac{x_1\otimes y_2- y_1\otimes x_2}{2}) \,.
\end{align*}

\ssegment{110522ee}{Definition}
Similarly, we define the \emph{Heisenberg group extension} (or simply \emph{Heisenberg extension}), $H$, as follows. As a representation,
$
H:= E_1 \oplus E_2 \oplus E_3 \;.
$
We give $H$ the structure of a group, unipotent of depth two, by setting
$$
X\cdot Y:= (X_1+Y_1, X_2+Y_2, X_3 +Y_3 + X_1 \otimes Y_2)\;.
$$

\ssegment{110522f}{}
If we write, for $x\in \hh$,
$$
e^x := (x_1, x_2, x_3+\frac{x_1 \otimes x_2}{2}) \;,
$$
then a straightforward verification shows that
$
x \mapsto e^x
$
defines an equivariant isomorphism
$
(\hh, \star) \xrightarrow{\cong} (H, \cdot) \;.
$

%%%%%%%%%%%%%%%%%%%%%%%%%%%%%%%%%%
\segment{1110182}{Coboundary maps in nonabelian group cohomology}%%%%%%%%%%%%%%%%%%%%
%%%%%%%%%%%%%%%%%%%%%%%%%%%%%%%%%

\ssegment{110523a}{}
Suppose $G$ acts on a group $U$. Then we define pointed sets $\m C^i, \m Z^i, \m B^i, \m H^i$ for small $i$ and maps $d^0:\m C^0 \to \m C^1$, $d^1:\m C^1 \to \m C^2$ as follows. $\m C^i$ is the set of continuous functions $G^i \to U$ together with the distinguished point given by $\si \mapsto e$, which we denote again by $e$. We define $d^0$ by
$$
(du)(\si) = u\si(u^{-1}) \;;
$$
and $d^1$ by
$$
dc(\si_1, \si_2) = c(\si_1\si_2)(\si_1c(\si_2))^{-1}c(\si_1)^{-1} \;.
$$
This allows us to define $\m H^0:= (d^0)^{-1}(e) = U^G$, and $\m Z^1:= (d^1)^{-1}(e)$. We define an action of $U$ on $\m Z^1$ by 
$$
(uc)(\si)=uc(\si)\si(u^{-1}) \;.
$$
Finally, $\m H^1$ is the set of orbits for this action.

\ssegment{110523b}{}
Now let 
$$
1 \to V \xrightarrow{i} U \xrightarrow{p} W \to 1
$$
be a short exact sequence of groups with $G$ action, suppose $V$ is central in $U$, and, supposing $p$ admits a (continuous) section, we fix one arbitrarily and call it  $q$. Maps
\begin{align*}
\delta:\m H^0(W) \to \m H^1(V) && \text{ and } && \Delta:\m H^1(W) \to \m H^2(V)
\end{align*}
may be defined as follows. Given $w \in \m H^0(W)$, lift $w$ to a $u \in U$ and then apply $d$:
$
\delta(w) := d(u) \;.
$
Similarly, given $c \in \m Z^1(W)$ use $q$ to lift $c$ to a cochain $b \in \m C^1(U)$ and then apply $d$:
$
\Delta(c):=d(b) \;.
$

\ssegment{110523c}{Proposition}
$\delta$ and $\Delta$ are well defined, independently of $q$, and give rise to an exact sequence of pointed sets
\[
e \to V^G \to U^G \to W^G \to \m H^1(V) \to \m H^1(U) \to \m H^1(W) \to \m H^2(V)
\,.
\]

\medskip
\noindent
Our conventions for the boundary maps follow \cite{kimii}. With different (but equivalent) definitions, this is proved in \cite{corloc}.

%%%%%%%%%%%%%%%%%%%%%%%%%
\segment{111213b}{Relationship to cup product}%%%%%
%%%%%%%%%%%%%%%%%%%%%%%%%

\ssegment{110523d}{}
Returning to the situation and the notation of \ref{110522b}, we recall (for instance, from \cite{Tate}), that the cup products
\begin{align*}
\m H^1(E_1) \times \m H^1(E_2) \to \m H^2(E_3) \;,
&&
\text{and}
&&
\m H^1(E_2) \times \m H^1(E_1) \to \m H^2(E_3)
\end{align*}
are given on the level of cocycles by the formulas
$$
(\ka_1 \cup \ka_2)(\si, \tau) = \ka_1(\si) \otimes \si\ka_2(\tau) \;,
$$
and
$$
(\ka_2 \cup \ka_1)(\si, \tau) = \si \ka_1(\tau) \otimes \ka_2(\si) \;.
$$
Thus defined, we have
$
\ka_1 \cup \ka_2 = - \ka_2 \cup \ka_1 \;.
$
Indeed, if we define $\ka \in \m C^1(E_3)$ by
$
\ka(\si)= \ka_1(\si) \otimes \ka_2(\si)
$
then
$$
\ka_1 \cup \ka_2 + \ka_2 \cup \ka_1 = d\ka
$$
in $\m Z^2(E_3)$.

\ssegment{110522g}{Proposition}
Let $G$ be a group, and $E_1, E_2$ representations over a ring containing $1/2$. Let $E_3=E_1\otimes E_2$, let $H$ denote the Heisenberg extension (\ref{110522ee}), and let
$$
\Delta: \m H^1(G, E_1 \oplus E_2) \to \m H^2(G, E_3)
$$
denote the second coboundary map (\ref{110523b}) associated to the exact sequence
$$
1 \to E_3 \to H \to E_1 \oplus E_2 \to 1 \;.
$$
Then $\Delta$ is related to the cup product by the formula
$
\Delta(\ka_1, \ka_2) = \ka_2 \cup \ka_1 \;.
$

\begin{proof}
This is merely a verification. Let
$
q: H \from E_1 \oplus E_2
$
denote the structural section. We note that while $q$ is not a homomorphism for the multiplicative structure of $H$, it is a homomorphism for the additive structure of $\hh$, which we now substitute for $H$. Fix $ (\ka_1, \ka_2) \in \m Z^1(G, E_1 \oplus E_2) $ and let $ \ka \in \m C^1(G, \hh) $ denote the lift of $(\ka_1, \ka_2)$ given by $q$. We then have $ \Delta(\ka_1, \ka_2) = d\ka \;, $ so for $\si, \tau \in G$, we have
\begin{align*}
\Delta(\ka_1, \ka_2)(\si, \tau) 	&= (\ka(\si)+\si\ka(\tau)) \star (-\si\ka(\tau)) \star (-\ka(\si)) \\
								&= -\frac{[\ka(\si),\si\ka(\tau)]}{2} \;,
\end{align*}
which has preimage
$$
\frac{\ka_2 \cup \ka_1 - \ka_1 \cup \ka_2}{2}(\si, \tau) = (\ka_2 \cup \ka_2) (\si, \tau)
$$
in $E_3$, indeed.
\end{proof}

%%%%%%%%%%%%%%%%%%%%%%%%%%%%%%%%%%
\segment{111213c}{Interpretation in terms of extensions}%%%%%%%
%%%%%%%%%%%%%%%%%%%%%%%%%%%%%%%%%

\ssegment{110525a}{}
Let $G$ be a group, and $\rho_i: G \to \GL(L_i)$, $i=1, \dots n$, representations over a ring $A$. Let
$
L_{i,j}= \Hom_{\Mod(A)}(L_j, L_i) \;,
$
and let
$
L_{i,j}^*= \Isom_{\Mod(A)}(L_j, L_i) \;.
$
Let $U$, $T$, $P$ denote the groups of block matrices:
\[ U=
\begin{pmatrix}
I & L_{1,2} & L_{1,3} & \cdots \\
& I & L_{2,3} & \cdots \\
&& I \\
&&& \ddots
\end{pmatrix}
\;,
\hspace{5mm}
T=
\begin{pmatrix}
L_{1,1}^* & & &  \\
& L_{2,2}^* & &  \\
&& L_{3,3}^* \\
&&& \ddots
\end{pmatrix}
\;,
\]
\[
P=
\begin{pmatrix}
L_{1,1}^* & L_{1,2} & L_{1,3} & \cdots \\
& L_{2,2}^* & L_{2,3} & \cdots \\
&& L_{3,3}^* \\
&&& \ddots
\end{pmatrix} \;.
\]
Define
$
\rho: G \to T
$
by
$$
\rho(\si) =
\begin{pmatrix}
\rho_1(\si) \\
& \rho_2(\si) \\
&& \ddots
\end{pmatrix}
\;.
$$
We let $G$ act on $U$ componentwise. This action is compatible with the group structure on $U$. Indeed, given $\si \in G$ and $x \in U$, we have
$$
\si x = \rho(\si)x\rho(\si)\inv \;.
$$
Let $\Hom^\rho(G, P)$ denote the set of homomorphisms whose composite with the projection $P \twoheadrightarrow T$ is equal to $\rho$. Given $\ka \in \m{Z}^1(G, U)$ and $\si \in G$, we let
$$
\rho_\ka(\si) = \ka(\si) \rho(\si) \;.
$$
Then the map
$
\ka \mapsto \rho_\ka
$
defines a bijection
$$
\m{Z}^1(G, U) \xrightarrow{\cong} \Hom^\rho(G, P) \;.
$$
The action of $U$ on the left corresponds to the action on the right by conjugation. The quotient $[\Hom^\rho(G,P)/U]$ parametrizes representations $E$ of $G$ equipped with a filtration
 $\Fil$ by subrepresentations, and an isomorphism
$$
\gr_{\Fil}E \xrightarrow{\cong} \bigoplus L_i \;.
$$
We call such objects \emph{$\rho$-framed representations}. We've thus constructed a bijection as follows.

\subsection*{Proposition}
$ \m{H}^1(G, U) = \Hom^\rho(G,P)/U = $ the set of isomorphism classes of $\rho$-framed representations of $G$.

\ssegment{110525b}{}
The group $U$ of paragraph \ref{110525a} specializes in the case $n=3$ to the Heisenberg extension
$$
1 \to L_{1, 3} \xrightarrow{i} H \xrightarrow{p} L_{1,2} \oplus L_{2,3} \to 1 \;.
$$
We have
$
\m{H}^i(G, L_{i,j}) = \Ext^i_G(L_j, L_i) \;,
$
so the long exact sequence of \ref{110523c} gives rise to an exact sequence
$$
\m{H}^1(G, H) \xrightarrow{p_*} \Ext^1(L_2, L_1) \oplus \Ext^1(L_3, L_2) \xrightarrow{\Delta} \Ext^2(L_3, L_1) \;.
$$
Given an element $(E', E'')$ of the middle term, paragraph \ref{110525a} shows that we can interpret $(p_*)\inv(E', E'')$ as the set of isomorphism classes of diagrams like the diagram on the left.
\[\xymatrix{
& 0 \ar[d] & 0 \ar[d] &&& 0 \ar[d]  \\
& L_1 \ar[d] \ar@{=}[r] & L_1 \ar[d] &&& L_1 \ar[d] \\
0 \ar[r] & E' \ar[r] \ar[d] & E \ar[d] \ar[r] & L_3 \ar[r] \ar@{=}[d] & 0 & E' \ar[dr] \\
0 \ar[r] & L_2 \ar[r] \ar[d] & E'' \ar[r] \ar[d] & L_3 \ar[r] & 0 & & E'' \ar[r] & L_3 \ar[r] & 0\\
& 0 & 0 
}
\]
On the other hand, $\Delta(E', E'')$ corresponds to the two-step extension on the right.
%\mar{proofs?}
%\[\xymatrix{
%0 \ar[d] \\
%L_1 \ar[d] \\
%E' \ar[dr] \\
%& E'' \ar[r] & L_3 \ar[r] & 0
%}
%\]

%%%%%%%%%%%%%%%%%%%%%%%%%%%%%%%
\section{Review of mixed Tate duality} \label{111018}%%%%%%%
%%%%%%%%%%%%%%%%%%%%%%%%%%%%%%%

\segment{111010a}{Definition}
A \emph{weakly mixed Tate category} over a field $k$ consists of a Tannakian category $T$ over $k$, plus an object $k(1)$, such that the $k(i):= k(1)^{\otimes i}$ are pairwise nonisomorphic objects of rank one, and every simple object is isomorphic to $k(i)$ for some $i \in \ZZ$. Recall from paragraph 2.1 of Beilinson-Deligne \cite{beilinson-deligne} that a \emph{mixed Tate category} is a weakly mixed Tate category satisfying
$
\Ext^1(k(0), k(n)) =0
$
whenever $n \le 0$. Given weakly mixed Tate categories $T$, $T'$, and a morphism of Tannakian categories $f:T \to T'$, we say that $f$ is \emph{mixed Tate} if $f(k(1)) = k(1)$.

\segment{111026}{}
Every object of a mixed Tate category is equipped with a canonical finite increasing filtration, which, following well motivated conventions, is denoted by $W$ and indexed by the even integers. If $E$ is an object, then each $\gr_{2n}^W E$ is isomorphic to a sum of $k(-n)$'s.

Conversely, if $T$ is merely a Tannakian category, and $k(1)$ is an object such that the $k(i):= k(1)^{\otimes i}$ are distinct simple objects, then the full subcategory $\langle k(1) \rangle_\m{mt}$ whose objects are those objects of $T$ which admit an increasing filtration $W$ indexed by the even integers, such that each $\gr^W_{2i}$ is isomorphic to a sum of $k(-i)$'s, is a mixed Tate category.

\segment{111010b}{}%%%%%%%%%
A mixed Tate category has a canonically defined fiber functor $\om$ over $k$, given by
\[
\om(E) = \bigoplus_{n \in \ZZ} \Hom(k(-n), \gr ^W_{2n}E)
\,.
\]
Its grading, given by
\[
\om_n(E) = \Hom(k(-n), \gr^W_{2n}E)
\]
gives rise to a semidirect product decomposition
\[
B:= \underline \Aut^\otimes (\om) = U \rtimes \Gm
\]
with $U$ prounipotent. Moreover, if we let $L$ denote the Lie algebra of $U$ and we let $L = \bigoplus L_i$ denote its decomposition under the action of $\Gm$, then $L_i = 0 $ for $i \ge 0$. The resulting equivalence of categories $T \to \Rep B$ sends $k(i)$ to the character of $\Gm$ of weight $-i$. 
%This sign discrepancy is an annoyance, which we attempt to alleviate by introducing the notation
%\[
%k_i := k(-i)
%\,.
%\]

Let $A$ denote the coordinate ring of $U$. Then under the action of $\Gm$, $A$ too decomposes as a direct sum $A = \bigoplus A_i \,,$ with $A_0=k$, and $A_i = 0$ for $i < 0 $. An element of the graded piece $A_i$ may fruitfully be interpreted as an appropriate equivalence class of certain ``framed" objects. For simplicity, we restrict attention to the case $i=2$, focusing on those points that we will need to use below. For a more complete discussion (and a slightly different formulation) see segment 3.2 of Goncharov's \textit{Multiple polylogarithms and mixed Tate motives} \cite{GonMPMTM}.

\segment{111010c}{}
If $E$ is an object of $\Rep B$, we let
$
E = \bigoplus E_i
$
denote the decomposition of the underlying representation of $\Gm$ into isotypic components. In terms of the $E_i$'s, the weight filtration is given by 
\[
W_{2n} E := \bigoplus_{i \le n} E_i
\,,
\]
so that
\[
\gr^W_{2n}E = E_n
\]
is a representation of $B$ in which $U$ acts trivially, and $\Gm$ acts with weight $n$.

\segment{111018c}{}
Let $F=(F_i)_{i\in I}$ be a family of objects of $T$ indexed by a subset $I \subset 2\ZZ$, and let $E$ be an object.

\subsection*{Definition}
An \emph{$F$-frame on $E$} is an isomorphism
\[
\gr^W_{n}E = F_n
\]
for each $n \in I$. We call the data consisting of $E$ together with an $F$-frame, \emph{an $F$-framed object}, or simply \emph{a framed object}. When $E$ is already in the picture and we wish to emphasize $F$, we say that $E$ \emph{is framed by $F$}. Given $F$-frames on objects $E$, $E'$, we say that a morphism $\phi: E \to E'$ is \emph{framed} if the associated graded pieces $\gr^W_n \phi$ for $n \in I$, are compatible with the isomorphisms defining the frames.

\medskip
\label{111010d}
%\mar{111010d}
\noindent
We let $(k(0),*,k(-2))$ denote the family of objects indexed by $2\ZZ \setminus \{2\}$
%\mar{The set of even integers except 2}
 given by 
\[
(k(0),*,k(-2))_i = 
\left\{
\begin{matrix}
k(0) & \mbox{ for } i= 0 \\
k(-2) & \mbox { for } i= 4 \\
0 & \mbox{ otherwise.}
\end{matrix}
\right.
\]

\segment{111010e}{}
We construct an abelian group, which we denote by $\Aa_2$. Underlying $\Aa_2$ is a partition of the set of $(k(0),* , k(-2))$-framed objects: two $(k(0),* , k(-2))$-framed objects $E$, $E'$ are to be considered equivalent whenever there exists a framed morphism $E \to E'$. For the group structure we let the semisimple object $k(0)\oplus k(-2)$ serve as origin. Let $\nabla$ denote the counterdiagonal in $k(-2)^{\oplus 2}$ (= kernel of the sum of the two projections). Given representatives $E$, $E'$ of elements of $\Aa_2$, let $\til E$ be the kernel of the natural surjection
\[
E \oplus E' \surj k(-2)^{\oplus 2} / \nabla
\,.
\]
We let $K$ denote the kernel of the sum
\[
W_0 \til E = W_0 E \oplus W_0 E' = k(0) \oplus k(0) \surj k(0)
\,
\]
of the two projection. Finally, to define addition, we set
\[
[E] + [E'] = [\til E / K]
\,.
\]
%\mar{what's the frame?}

\segment{111011}{}
%$\Aa_2$ fits naturally into a triangle of abelian groups like so,
%\[
%\xymatrix{
%\Aa_2 \ar[rr]^-\tau \ar[dr] _-\nu
%&&
%A_2 \ar[dl]^-\mu
%\\
%&
%A_1 \otimes A_1
%}
%\]
Consider an object $E$ representing an element of $\Aa_2$. Writing $k_i$ for the ground field $k$ regarded as vector space in graded degree $i$, $E$ decomposes under the $\Gm$ action as
\[
E = k_0 \oplus E_1 \oplus k_2
\,.
\]
The coaction map
$
A \otimes E \from  E
$
decomposes under the action of $\Gm$ as a direct sum of six maps, of which three are uniquely determined; the remaining three we denote by $a$, $b$, $c$ as follows:
\begin{align*}
A_1 \otimes E_{1} \xfrom {b} k_2
&&
A_1 \otimes k_{0} \xfrom {a} E_{1}
\end{align*}
\[
A_2 \otimes k_{0} \xfrom{c} k_2
\,.
\]
These maps are subject only to the condition that the square
\[
\xymatrix{
%1
A_1 \otimes A_1 														&
A_2 					\ar[l]_-\nu 									\\
%2
A_1 \otimes E_{1} 	\ar[u]^ {1_{A_1} \otimes a} 					&
k_2 						\ar[u]_-{c} \ar[l]^-{b}			}
\]
in which $\nu$ denotes comultiplication, commutes. Stated more geometrically, $E$ gives rise to elements $a \in A_1 \otimes E_{1}^\lor$, $b\in A_1 \otimes E_{1}$, and $c \in A_2$. If $C$ is a $k$-algebra and $x \in U(C)$, then the action of $x$ on $C \otimes E$ is given by the block matrix
\[
\rho(x) = 
\left(
\begin{matrix}
1 &a(x)&c(x)\\
&1_{E_{1}}&b(x)\\
&&1
\end{matrix}
\right)
\,.
\]
Matrix multiplication prescribes the conditions
\begin{align*}
a(xy) &=a(x)+a(y) \\
b(xy) &=b(x)+b(y) \\
c(xy) &=c(x)+c(y)+a(x)b(y)
\end{align*}
on $a$, $b$, $c$, in terms of an arbitrary pair of points $x,y \in U(C)$.

\segment{111017c}{}
Framed objects and framed morphisms as in segment \ref{111010e}, form a category, as do data of the form $(E_{1}, a, b, c)$ as in segment \ref{111011}, and the construction of segment \ref{111011} gives rise to an isomorphism of categories.

\segment{110316}{Claim}
In the notation of segments \ref{111010e}--\ref{111011}, the formula
$
[E]  \longmapsto c 
$
defines an isomorphism
$
\Aa_2  \liso  A_2
\,.
$

\medskip
\noindent
The remainder of segment \ref{110316} is devoted to verification of this statement. We omit the verification of additivity.

\ssegment{111017a}{}
It is clear that the map is well defined on equivalence classes: recall that a morphism of representative objects
$
E \to E'
$
is required to be \emph{framed}; through the frame, such a morphism gives rise to a commuting square
\[
\xymatrix{
A_2 \ar@{=}[r] & A_2 \\
k_2 \ar@{=}[r]  \ar[u] & k_2 \ar[u]
}
\]
from which it follows that $c'=c$.

\ssegment{111013}{}\textit{Surjectivity.}
Fix arbitrarily an element $c \in A_2$ and write
\[
\nu(c) = \sum_{i=1}^n b_i \otimes a_i
\,.
\]
%Before we define a framed object $E$ whose class in $\Aa_2$ maps to $c$, we note that $T$ (\ref{}) is canonically isomorphic to the category of finite dimensional graded comodules over $A$. Through this isomorphism, $k(i)$ becomes $k$ in degree $i$ with trivial coaction. 
We define $E$ by setting $E_{1}:= k^n$ with $i^\m{th}$ standard basis vector denoted by $e_i$, and we define the three variable components of the coaction map $a$, $b$, $c$ by
\begin{align*}
\sum b_i \otimes e_i \underset b \longmapsfrom 1
\,,
&&
a_i \underset a  \longmapsfrom e_i
\,,
&&
\mbox{and}
&&
c \longmapsfrom 1
\,,
&&
\mbox{respectively.}
\end{align*}

\ssegment{111017b}{} \textit{Injectivity.}
Suppose $c=0$. Define a framed object $E'$ as follows. Let $E'_{1}$ be the kernel of $a$:
\[
0 \to E'_{1} \to E_{1} \to A_1
\,.
\]
Set $c' =0$ and $a'=0$. To define $b'$, note that the composite
\[
\xymatrix{
A_1 \otimes A_1 \\
A_1 \otimes E_{1} \ar[u]^{1_{A_1}\otimes a} & k \ar[l]^-b
}
\]
is zero, so that $b$ factors through $A_1 \otimes E'_{1}$. Finally, the maps
\[
E_{1} \hookleftarrow E'_{1} \to 0
\]
of vector spaces, define framed morphisms of framed objects. Together, these witness the desired equivalence
$
E \sim 0
$
in $\Aa_2$.

\segment{111018d}{}
Similarly, but more simply, we may construct an abelian group $\Aa_1$ as follows. We consider the set of $(k(0), k(-1))$-framed objects, we stipulate that framed objects $E$, $E'$ are \emph{equivalent} whenever there exists a framed morphism $E \to E'$, and we let $\Aa_1$ be the set of equivalence classes. But an object \emph{framed by $(k(0),k(-1))$} is no different from an extension
\[
0 \to k(0) \to E \to k(-1) \to 0 
\,,
\]
and a \emph{framed morphism} is the same as an isomorphism of extensions. So
\[
\Aa_1 = \Ext^1_T(k(0),k(1)) = H^1(T, k(1))
\]
tautologically.

\segment{111018e}{}
In the situation and the notation of segment \ref{111018d}, the coaction is only a map
\[
A_1 \otimes k_0 \from k_1
\,.
\]
If we let $a$ denote the corresponding element of $A_1$, then the formula
$
E \longmapsto a
$
defines an isomorphism
$
\Aa_1 \liso A_1
\,.
$

\segment{120103a}{}%%%%%%%%%%
We now discuss how to translate the product $\mu :A_1 \otimes A_1 \to A_2$ into the language of framed objects. If $E'$, $E''$ are representations of $B$, then we have an isomorphism of $B$-representations
\[
\gr_n^W(E' \otimes E'') = \bigoplus_{i+j=n} (\gr^W_i E')\otimes (\gr_j^W E'')
\,.
\]
(Here all indices range over the even integers. To see this we need only recall that $U$ acts trivially on $W$-graded pieces.) Thus, if $E'$, $E''$ are $(k(0),k(-1))$-framed objects, $E:= E' \otimes E''$ naturally has the structure of a $(k(0), *, k(-2))$-framed object. We denote it by $E = E' \varoast E'' $. The coactions of $A$ on $E'$ and $E''$ are determined by elements $a', a'' \in A_1$. In terms of $a', a''$, the coaction of $A$ on $E$ is given by the following square.
\[
\xymatrix
@ C=2cm
{
%1
A_1 \otimes A_1																&
A_2					\ar[l]_-\nu												\\
%2
A_1 \otimes k^2		\ar[u]^-{1_{A_1} \otimes (a'' \oplus a')}					&
k						\ar[l]^-{a' \otimes e_1 + a'' \otimes e_2} \ar[u]_{a'a''}		}
\]
This establishes the following two propositions.

\ssegment{120103b}{Proposition}%pppppppp
In terms of framed objects, the multiplication map $\mu: A_1 \otimes A_1 \to A_2$ is given by
\[
[E'] \otimes  [E''] \to [E' \varoast E'']
\,.
\]

\ssegment{120103c}{Proposition}%pppppppp
The composite
\[
A_1 \otimes A_1 \xto \mu A_2 \xto \nu A_1 \otimes A_1
\]
is given by
\[
a' \otimes a'' \mapsto a' \otimes a'' + a'' \otimes a'
\,.
\]

%%%%%%%%%%%%%%%%%%%%%%%%%%%%%%%%%%%
\section{Heisenberg torsors in mixed Tate categories}%%%%%%%
\label{120201-1}%%%%%%%%%%%%%%%%%%%%%%%%%
%\mar{120201-1}%%%%%%%%%%%%%%%%%%%%%%%%%%
%%%%%%%%%%%%%%%%%%%%%%%%%%%%%%%%%%%

We continue to work with a mixed Tate category $T$ over a field $k$, and we preserve all notations introduced in \S\ref{111018}. In particular, $B = U \rtimes \Gm$ denotes the canonical fundamental group.
%We often confuse an object with the corresponding representation; indeed, throughout this section we may assume $T = \Rep B$.

\segment{111018a}{Definition}%
We define the \emph{Heisenberg $T$-group $H_T$} to be the Heisenberg group extension in the sense of segment \ref{110522ee} with $E_1=E_2=k(1)$, and $E_3=k(2)$.  $H$ is a unipotent group object of $T$ in the sense of Example 5.14(iv) of Deligne \cite{Deligne89}, or, which is canonically the same thanks to the canonical fiber functor, the usual Heisenberg group over $k$, equipped with an action of $B$, or, which is the same, a graded action of $U$. Elaborating on this terminology, we call an object framed by $(k(n), k(n-1), k(n-2))$ (for any $n\in \ZZ$) a \emph{Heisenberg $T$-torsor}.

\segment{111018f}{}%%%%%%%%%%%%
Using a fiber functor, the first cohomology set $\m H^1(T,G)$ of a group $T$-object $G$, may be defined in terms of explicit cocycles as in segment \ref{1110182}. So in particular, using the canonical fiber functor, we have
\[
\m H^1(T, H_T) = \m H^1(B, H_T)
\,.
\]
According to segment \ref{110525a}, the result is nevertheless independent of the choice of fiber functor; moreover, $\m H^1(T, H_T)$ may be identified with the set of isomorphism classes of Heisenberg $T$-torsors.  Our goal in this section is to construct a commutative diagram with exact columns, like so.
\[
\xymatrix{
\Ext^1_T(k(0),k(2))			\ar[d]_-g 						&
\m H^1(T, k(2)) 				\ar@{=}[l] \ar[d] 				\\
A_2							\ar[d]_-\nu						&
\m H^1(T,H_T)					\ar[l]_-{\la_2}  \ar[d]			\\
A_1 \otimes A_1				\ar[d]_-D 						&
\m H^1(T, k(1) \oplus k(1)) 	\ar[l]_-{\la_1}  \ar[d]^-\Delta	\\
\Ext^2_T(k(0),k(2))											&
\m H^2(T, k(2))					\ar@{=}[l]						}
\]
The right column is simply the long exact sequence of nonabelian cohomology (\ref{110523c}) applied to the short exact sequence
\begin{iquation}
\label{111019}
0 \to k(2) \to H_T \to k(1) \oplus k(1) \to 1
\,.
\end{iquation}
Segments \ref{111018d} and \ref{111018e} together compose an isomorphism
\[
A_1 = \Aa_1 = \m H^1(T, k(1))
\]
from which we obtain $\la_1$ as the universal bilinear map. Similar, but somewhat more complicated, is the construction of $\la_2$. According to segment \ref{110525a}, an element of $\m H^1(T, H_T)$ may be represented by an object $E$, framed by $(k(0), k(-1), k(-2))$ (\ref{111018c}). We define its image in $A_2$ through the isomorphism $A_2 = \Aa_2$ of segment \ref{110316}, by forgetting the isomorphism
\[
\gr^W_{2}E = k(-1)
\,.
\]
In addition to $\la_2$, this gives us $g$. To define $D$, we recall from Proposition \ref{110522g} that $\Delta$ is bilinear, hence factors through $A_1 \otimes A_1$.

\segment{111019a}{}%%%%%%%%%%%%%%%
Verification of the commutativity of the central square is straightforward, and we omit it. Commutativity of the remaining two squares is clear. Exactness of the left column at $A_2$ goes as follows: given $c \in \ker \nu$, we endow $k(0)\oplus k(2)$ with the coaction given by $c$ to obtain an extension of $k(0)$ by $k(2)$.

\segment{111022a}{}%%%%%%%%%%%%%
 It remains to prove that the left column is exact at $A_1 \otimes A_1$. Given $d \in A_1 \otimes A_1$, there exist a vector space $E_{1}$ and elements $a \in A_1 \otimes E_{1}^\lor$, $b \in A_1 \otimes E_{1}$, such that $ab=d$. Indeed, these may be constructed out of the choice of a decomposition $d = \sum_1^n b_i \otimes a_i$ by setting $E_{1} = k^n$,
\begin{align*}
a =
\begin{pmatrix}
a_1, \dots, a_n
\end{pmatrix}
\,,
&&
\text{and}
&&
b =
\begin{pmatrix}
b_1 \\
\vdots \\
b_n
\end{pmatrix}
&&
.
\end{align*}
Define objects $E'$, $E''$ with $E'_{1} = E''_{1} := E_{1} $ by creating the graded vector spaces
%\mar{multiple small problems here}
\begin{align*}
E' = k_0 \oplus E_{1}
&&
E'' = E_{1} \oplus k_2
\,,
\end{align*}
and endowing $E'$, $E''$ with the graded $A$ coaction given by $a$, $b$, respectively. Then the image of $d$ in $\Ext^2$ may be represented by the two-step extension
\[
0 \to k(0) \to E' \to E'' \to k(-2) \to 0 \, .
\]
It is nullhomotopic if and only if there exists a commutative diagram with exact rows like so.
\[
\xymatrix{
& 0 \ar[d] & 0 \ar[d] \\
& k(0) \ar[d] \ar@{=}[r] & k(0) \ar[d] \\
0 \ar[r] & E' \ar[r] \ar[d] & E \ar[d] \ar[r] & k(-2) \ar[r] \ar@{=}[d] & 0 \\
0 \ar[r] & E_{1} \ar[r] \ar[d] & E'' \ar[r] \ar[d] & k(-2) \ar[r] & 0 \\
& 0 & 0
}
\]
This makes $E$ into a $(k(0), * ,k(-2) )$-framed object, and, through the isomorphism \ref{110316}, gives us an element of $A_2$ which lies in the preimage of $d$.

%%%%%%%%%%%%%%%%%%%%%%%%%%%%%%%%%
\section{The polylogarithmic quotient}\label{120129-1}%%%%%
%\mar{120129-1}%%%%%%%%%%%%%%%%%%%%%%%%%
%%%%%%%%%%%%%%%%%%%%%%%%%%%%%%%%%

\segment{0}{}
We now turn our attention to the thrice punctured line $X:=\PP^1_{\QQ} \setminus \{0,1,\infty\}$. For $x$ either a $\QQ$-valued point or a tangent vector at one of the punctures, we denote by $\pioneinfty(X,x)$ the unipotent fundamental group of $X$ at $x$ as in Deligne-Goncharov \cite{DelGon} --- a unipotent group object of the Tannakian category of mixed Tate motives over $\ZZ$ with $\QQ$ coefficients. For $n\in\NN$ let $Z^n\subset\pioneinfty(X,x)$ denote the $n^\m{th}$ step in the
  descending central filtration and $\pi_1^{[n]}(X,x):=\pioneinfty(X,x)/Z^{n+1}$ the
  quotient. The main goal of this section is to show that for $n=2$ we have
\begin{align*} \label{eq1.1}
     \pi_1^{[2]}(X,\tbp) \; = \; 
       \begin{pmatrix}
          1 & \QQ(1) & \QQ(2) \\
          0 &   1    & \QQ(1) \\
          0 &   0    &   1  
       \end{pmatrix}.
\end{align*}
This isomorphism is essentially obtained from the monodromy of the 
dilogarithm. Although it is clear from the Betti or de Rham realization that
we have a short exact sequence
\[
     1\to \QQ(2) \to\pi_1^{[2]}(X,\tbp) \to \QQ(1)^2\to 1
\]
which, at the level of algebraic groups, is the Heisenberg extension, it is
less obvious that this sequence splits at the level of motives (i.e.\ if we
consider the associated sequence of Lie algebras and forget the algebra
structure).  To prove this, we construct a faithful unipotent representation
of $\pi_1^{[2]}(X,\tbp)$ on a motive of the form
\[
    V_2= \QQ(1)\oplus\QQ(2)\oplus\QQ(3).
\]

\segment{0}{}%%%%%%%%%%%%%
Let $x$ be an arbitrary base point. Following Deligne (\cite{Deligne89}, \S
16.11), we define a certain quotient $U(x)$ of $\pioneinfty(X,x)$ as follows. The
inclusion $X\inj \GG_m$ defines a morphism
\[
\lambda_x:\pioneinfty(X,x) \to \pioneinfty(\GG_m,x)=\QQ(1)
\,.
\]
Let $N_x$ be the kernel of
$\lambda_x$ and $N_x'\subset N_x$ the derived subgroup. We define
$U(x):=\pioneinfty(X,x)/N_x'$. Clearly, we obtain a short exact sequence
\begin{iquation} \label{eq1.2}
  1 \to U_1(x) \to U(x) \to \QQ(1) \to 1,
\end{iquation}
where the kernel $U_1(x):=N_x/N_x'$ is abelian. 

For $x=\vec{01}$ we simply
write $U:=U(\vec{01})$ and $U_1$, for $x=\vec{10}$ we write $U'$ and $U_1'$. 
The tangential base points yield canonical `monodromy' morphisms
\[
   \mu_0:\QQ(1)\to U, \qquad \mu_1:\QQ(1)\to U'.
\]
The first of these morphisms defines a splitting of \eqref{eq1.2}, so we have
$U\cong U_1\rtimes\QQ(1)$. The image of the second actually lies in $U_1'$, so
we have a morphism $\mu_1:\QQ(1)\to U_1'$. 

We claim that there exists a canonical isomorphism $U_1\cong U_1'$. To see
this, we first note that $U'$ is isomorphic to the twist of $U$ by the
$U$-torsor $P:=U(\vec{01},\vec{10})$.
%\marginpar{Definition of $P$?} 
If we
restrict this isomorphism to $U_1$ and use the fact that $U_1$ is abelian, we
see that $U_1'$, as a twist of $U_1$, depends only on the image of $P$ on
$\GG_m$, i.e.\ on the $\QQ(1)$-torsor
$\bar{P}=\pioneinfty(\GG_m,\vec{01},\vec{10})$. But $\bar{P}$ is trivial and hence
$U_1'\cong U_1$. Moreover, since $\bar{P}$ has a unique trivialization, our
argument shows that this isomorphism is canonical. This proves the claim.

By the claim we have a morphism $\mu_1:\QQ(1)\to U_1$. For $n\geq 1$ we define
\[
    \nu_n:=(\ad\,\mu_0)^{n-1}(\mu_1):\QQ(n)\to U_1.
\]
Clearly, $\nu_1=\mu_1$. According to Proposition 16.13 of Deligne \cite{Deligne89}, $(\nu_1,\nu_2,\ldots)$ defines an
isomorphism 
\begin{iquation} \label{eq1.3}
    \Lie U_1 \cong \prod_{n\geq 1} \QQ(n).
\end{iquation}
Let us recall how this works for the de Rham realization. According to Segment 15.52 of Deligne \cite{Deligne89}, 
 $L:=\Lie\pi_1^{\rm dR}(X)$ is the nilpotent completion of the free Lie
algebra generated by $X_0,X_1$, where $(X_0,X_1)$ is the dual basis of
$(\omega_0=dz/z,\omega_1=dz/(z-1))$ of $H^1_{\rm dR}(X)$. One checks that
$X_0=\mu_0^\dR(1)$ and $X_1=\mu_1^\dR(1)$. From there it follows easily that
$\Lie N^\dR\subset L$ is the Lie ideal generated by $X_1$. Moreover, if
$e_0,e_1\in\Lie U^\dR$ denote the images of $X_0,X_1$, then
\[
        \Lie U^\dR = \gen{e_0;(\ad\,e_0)^{n-1}(e_1), n\geq 1}^{\rm top},
\]
and the elements $(\ad\,e_0)^{n-1}(e_1)$ commute with each other, completing the argument.

Now let us fix $n\geq 2$. We define $V_n\subset \Lie U^{[n+1]}$ as the
ideal generated by $e_0$ and $e_2:=[e_0,e_1]$. Then
$(e_0,e_2,\ldots,e_{n+1})$, with $e_k:=(\ad\,e_0)^{k-1}(e_1)$, is a basis of
$V_n$ and provides via \eqref{eq1.3} a decomposition
\begin{align*} \label{eq1.4}
    V_n\cong \oplus_{k=1}^{n+1}\,\QQ(k).
\end{align*}
 
Considering $V_n$ as a $\pi_1^{[n]}(X,\tbp)$-module, it corresponds to a
unipotent local system $\V_n$ on $X$ with $V_n=(\V_n)_{\tbp}$. For instance,
its de Rham realization is the trivial vector bundle
$\V_n^\dR=\Oo_X\otimes_\QQ V_n^\dR$ together with the unipotent connection
\[
    \nabla := 1\otimes d - X_0\otimes dz/z - X_1\otimes dz/(z-1)
\]
(here we consider $X_0,X_1$ as elements of $\End(V_n^\dR)$). In terms of the
basis $(e_0,e_2,\ldots,e_{n+1})$ this means that the connection matrix is
\[\begin{pmatrix}
       0     &   0      & \cdots & \cdots   &   0  \\
    \omega_1 &  0       & \cdots & \cdots   &   0  \\
       0     & \omega_0 & \ddots &          & \vdots \\
    \vdots   & \ddots   & \ddots & \ddots   & \vdots \\
       0     & \cdots   &   0    & \omega_0 &  0 \\
\end{pmatrix}.\]
One recognizes that $\V_n$ is the polylogarithmic local system described in
\cite{Hain92}.  

For $n=2$ we obtain a representation
\[
  \rho:\pi_1^{[2]}(X,\tbp) \to \begin{pmatrix}
          1 & \QQ(1) & \QQ(2) \\
          0 &   1    & \QQ(1) \\
          0 &   0    &   1  
       \end{pmatrix}
\]
which is easily seen to be an isomorphism. 

%%%%%%%%%%%%%%%%%%%%%%%%%%%%%%%%%%%%%%
\section{Explicit Chabauty-Kim theory for the multiplicative group}%%%
%%%%%%%%%%%%%%%%%%%%%%%%%%%%%%%%%%%%%%

%%%%%%%%%%%%%%%%%%%%%%%%%%%%%%%%%%%%
\segment{0}{The local unipotent $p$-adic Hodge morphism}%%%%
%%%%%%%%%%%%%%%%%%%%%%%%%%%%%%%%%%%%

\ssegment{110718a}{}%%%%%%%%%%%%%%%%%%
The unipotent $p$-adic Hodge morphism of the multiplicative group over $\Qp$ fits into a triangle like so.
\[
\xymatrix{
																		&
\Gm(\Zp)									\ar[dl]_-\ka \ar[dr]^-\al 		&
																		\\
\m H^1_\m f(G_\Qp; \pioneinftyet(\Gm)) 		\ar[rr]_-h					&
																		&
\pioneinftydr(\Gm)(\Qp)												}
\]
Here, $G_\Qp$ denotes the Galois group of $\overline\QQ_p$ over $\Qp$, $\pioneinfty(\Gm)$ denotes the unipotent fundamental group at $1$, $\ka$ denotes the unipotent \'etale Kummer map, and $\al$ denotes the unipotent Albanese map. Since $ \pioneinfty(\Gm) = \QQ(1) , $ the above triangle may be written as follows.
\[
\xymatrix{
																		&
\Zpstar										\ar[dl]_-\ka \ar[dr]^-\al 		&
																		\\
\m H^1_\m f(G_\Qp; \Qp(1))				\ar[rr]_-h					&
																		&
\Qp(1)_\textrm{dR}													}
\]
Our goal here is to analyze this triangle.

\ssegment{110718b}{}%%%%%%%%%%%%
We first compute $\m H^1(G_\Qp; \Qp(1))$. For each $n$ there's an isomorphism
\[
\QQ_p^*/(\QQ_p^*)^{p^n} \xrightarrow{\cong} \Zpn \times p\ZZ/p^{n+1}
\]
given by $ up^r \mapsto (r, \log u)$. (Here $\log$ denotes the $p$-adic logarithm, its values on units are well defined independently of choice of branch parameter; we need only declare the logarithm of $(p-1)^\m{st}$ roots of $1$ to be zero). This holds because $\QQ_p^* = \ZZ \times \ZZ_p^*$, because $\ZZ_p^*$ fits into a short exact sequence
\[
\xymatrix{
1 \ar[r] &
 (p) \ar[r]^-\exp &
  \ZZ_p^*  \ar[r] &
   \FF_p^* \ar[r] \ar@/_/[l] &
    1 
    }
\]
which is canonically split by the Teichmuller lift, and because for each $n$, $p^n$ is invertible in $\FF_p^*$. Subsequently, via the inverse system of Kummer exact sequences, we obtain $ \m H^1(G_\Qp; \Qp(1)) = \Qp \otimes_\Zp \underset{\from}\lim \QQ_p^*/(\QQ_p^*)^{p^n}  = \Qp \otimes_\Zp (\Zp \times p\Zp) = \Qp \oplus \Qp $.

\ssegment{110719a}{Proposition}%PPPPPPPPP
The isomorphism of paragraph \ref{110718b} restricts to
$$
\m H^1_\m f(G_\Qp; \Qp(1)) = 0 \oplus \Qp \;.
$$
Moreover, the unipotent $p$-adic Hodge morphism $h$ is given by this isomorphism.

\begin{proof}
This is explained in Remark II.4.4 of Berger \cite{Berger}.
\end{proof}

\ssegment{110719b}{}%%%%%%%%%%%
It follows from the commutativity of this triangle that $\al$ is given by the $p$-adic logarithm. We check this directly. The isomorphism
$$
\pioneinfty(\Gm; 1) \xrightarrow{\cong} \QQ(1)
$$
constructed in \cite{Deligne89}, \S14, has $p$-adic de Rham realization
$$
\pioneinftydr(\Gm;1) \xrightarrow{\cong} \Qp(1)_\textrm{dR}
$$
given by
$$
\ga \mapsto \left( \ga \left( \Oo^{(2)}, \left( 
\begin{matrix}
0 & 1\\
0 & 0
\end{matrix}
\right) \frac{dt}{t} \right) \right)_{1,2}
\;.
$$
On the other hand, for $z \in \Gm(\Zp)$ we have
$$
\al(z)\left( \Oo^{(2)}, \left( 
\begin{matrix}
0 & 1\\
0 & 0
\end{matrix}
\right) \frac{dt}{t} \right)
=
\left(
\begin{matrix}
1 & \int_1^z \frac{dt}{t} \\
0 & 1
\end{matrix}
\right)
\;.
$$
In Besser's formalism (\cite{Besser}), this is indeed the definition of the $p$-adic logarithm.

\segment{0}{The global Selmer variety}%%%%%%%%%%%%%%%%%%%%%%

\ssegment{110720a}{Proposition}%PPPPPPPPPPPPP
There is an isomorphism (constructed below)
$$
\m H^1(G_T, \Qp(1)) = \QQ_p^T \;.
$$
Moreover, the restriction map
$
\m H^1(G_T, \Qp(1)) \to \m H^1(G_\Qp, \Qp(1))
$
corresponds through this isomorphism and the isomorphism of paragraph \ref{110718b} to the map
\[
\QQ_p^T \to \Qp \oplus \Qp
\]
given by $(x_1, \dots, x_s, x) \mapsto	 (x, (\log q_1)x_1 + \cdots + (\log q_s)x_s)$, where $q_1, \dots, q_s$ are the elements of $S$.

\ssegment{110720b}{}%%%%%%%%%%%%
The construction goes as follows. Let $K_T$ denote the maximal unramified outside $T$ extension of $\QQ$, and let $\overline {\ZZ[T\inv]}$ denote the integral closure of $\ZZ[T\inv]$ in $K_T$.

\subsection*{Claim}
For each $n \in \NN$ the $T$-integral Kummer sequence
\[
1 \to \mupn \to \overline{\ZZ[T\inv]}^* \to \overline{\ZZ[T\inv]}^* \to 1
\]
is exact.

\begin{proof}
Only the surjectivity of the $(p^n)^\m{th}$ power map is in doubt. Let $a \in  \overline{\ZZ[T\inv]}^*$, let $v$ be a place of $K_T$ not above $T$, and let $w$ be a place of $K_T(\sqrt[p^n]{a})$ lying above $v$. Then the completion $(K_T)_v$ is isomorphic to the maximal unramified extension $\QQ_v^\m{ur}$ of $\QQ_v$, while the completion $(K_T(\sqrt[p^n]{a}))_w$ above it, is isomorphic to $\QQ_v^\m{ur}(\sqrt[p^n]{a})$. Since $p$ is coprime to $v$ and $a$ is a unit at $v$, the latter extension is unramified, hence trivial. Hence $\sqrt[p^n]{a} \in K_T$. Since moreover $\sqrt[p^n]{a}$ is integral, we have $\sqrt[p^n]{a} \in \overline{\ZZ[T\inv]}^*$, which establishes the claim.
\end{proof}

\medskip
%\ssegment{110720c}{}%%%%%%%%
\noindent
Meanwhile, we have 
\begin{align*}
\m H^1(G_T, \overline{\ZZ[T\inv]}^*) &= \m H^1(\Spec \ZZ[T\inv], \Gm) \\ &= \Pic \ZZ[T\inv] \\&=0 \,,
\end{align*}
 we have $\ZZ[T\inv]^* = \mu_2 \times \ZZ^T$, and we have
 \[
\m  H^1(G_T, \Qp(1)) = \Qp \otimes_\Zp \varprojlim \m H^1(G_T, \mupn)\,.
  \]
So when we apply $\m H^1(G_T, \cdot)$ to the $T$-integral Kummer sequence, we obtain the stated isomorphism.

The verification of our formula for the map is straightforward.

\ssegment{110720d}{Corollary}%%%%%%%%%%%%
There is an isomorphism
$
\m H^1_\m f(G_T, \Qp(1)) = \QQ_p^S \,,
$
such that the square
\[
\xymatrix{
\Gm(\ZSinv) \ar[r]^-\ka \ar[d] &
 \m H^1_\m f(G_T, \Qp(1)) \ar@{=}[d] \\
  \QQ^S \ar@{_(->}[r] &
  \QQ_p^S
}
\]
commutes. Moreover, the restriction map
$$
\m H^1_\m f(G_T, \Qp(1)) \to \m H^1_\m f(G_\Qp, \Qp(1))
$$
corresponds to the map
$
\AA^S_{\Qp} \to \AA^1_\Qp
$
given by 
\[
(x_1, \dots, x_s) \mapsto (\log q_1)x_1 + \cdots + (\log q_s)x_s \,. 
\]

%%%%%%%%%%%%%%%%%%%%%%%%%%%%%%%%%%%%%%
\section{Local computations}%%%%%%%%%
%%%%%%%%%%%%%%%%%%%%%%%%%%%%%%%%%%%%%%

\segment{110721a}{}%%%%%%%%%%%%%%%%%%%
We now turn our attention to the local picture, given in the notation ($X$, $\pioneet{n}$, $\pionedr{n}$) of the introduction, by the triangle
\[
\xymatrix{
	& X(\Zp) \ar[dl]_-\ka \ar[dr]^-\al 		&	\\
\m H^1_\m f(G_\Qp; \pioneet{2}) \ar[rr]_-h	&	& \pionedr{2}(\Qp)}
\]
in which $\ka$ now denotes the \textit{local} unipotent \'etale Kummer map.

\segment{110721b}{}%%%%%%%%%%%%%%%%%%%%%
See \cite{FurushoI}, \cite{FurushoII} for the definition of de Rham polylogarithms; and see Definition 2.32 of \cite{FurushoII} for the definition of \'etale plylogarithms which we will adopt. We recall that $\pioneinftydr$ has a canonical pair of generators giving rise to a canonical isomorphism
$
\pionedr{2} \xrightarrow{\cong} \UU_3
$
to the group of $3 \times 3$ upper triangular matrices with $1$'s on the diagonal.

\segment{110721c}{Proposition}%PPPPPPPPPP
The unipotent Albanese map in depth two (composed with the isomorphism of paragraph \ref{110721b}) is given by
$$
z \mapsto
\left(
\begin{matrix}
1 & \log z & -\Li z \\
0 & 1 & \log (1-z) \\
0 & 0 & 1
\end{matrix}
\right) \;.
$$
Here, $\log$ denotes the $p$-adic logarithm, and $\Li$ denotes the $p$-adic dilogarithm.

\begin{proof}%pppppppppp
Let $\UU(A,B)$ denote the free pro-unipotent group on the generators $A$, $B$, and let $\Qp<<A, B>>$ denote the ring of noncommutative power series in $A$, $B$. There is a natural embedding $\UU(A,B)(\Qp) \subset \Qp<<A,B>>$. An embedding
$$
\pioneinftydr(\Qp) \hookrightarrow \Qp<<A,B>>
$$
is defined in Notation 2.1 of \cite{FurushoII}. Its image is precisely $\UU(A,B)(\Qp)$. We let $\al_\infty$ denote the unipotent Albanese map. The logarithm depends on the choice of a branch parameter, for which we choose $0$.
%\mar{What's this branch parameter thing?}
According to Example 3.25 of \cite{FurushoI} and Theorem 2.3 of \cite{FurushoII}, the composite
$$
X(\Zp) \xto{\al_\infty} \pioneinftydr(\Qp) \xto{\sim} \UU(A,B)(\Qp) \subset \Qp<<A,B>>
$$
is given by
\begin{align*}
z \mapsto 1 &+ (\log z)A + \log (1-z) B + \frac{(\log z)^2}{2} A^2 \\
						& - \Li_2(z) AB + \{ Li_2(z) + (\log z)\log (1-z) \} BA + \frac{\{\log (1-z) \}^2}{2} B^2
\end{align*}
plus higher terms. Under the projection
$
\UU(A,B) \surj \UU_3
$
determined by
$$
A \mapsto
\left(
\begin{matrix}
1&1&0\\
0&1&0\\
0&0&1
\end{matrix}
\right)
\hspace{5mm}
\mbox{and}
\hspace{5mm}
B \mapsto
\left(
\begin{matrix}
1&0&0\\
0&1&1\\
0&0&1
\end{matrix}
\right)
\;,
$$
this element of $\UU(A,B)(\Qp)$ has image in $\UU_3(\Qp)$ as stated in the proposition, as one sees, for instance, by going through the corresponding Lie algebras.
\end{proof}

\segment{110722a}{}%PPPPPPP
We can give a similar formula for the unipotent \'etale Kummer map in terms of the \'etale dilogarithm.

\subsection*{Proposition}
The unipotent \'etale Kummer map in depth $2$:
$
X(\Zp) \to \m H^1_\m f(G_\Qp, H_\et)
\,,
$ 
is given by $ z \mapsto \ka_{2,z} $ where
$$
\ka_{2,z}(\si) = ( \ka_z(\si), \ka_{1-z}(\si), -\Li_z(\si)) \;.
$$

\medskip
\noindent
We note that this entails the cocycle condition on the \'etale dilogarithm:
$$
\Li_{2,z}(\si \tau) = \Li_{2,z}(\si) + \si(\Li_{2,z}(\tau)) - \ka_z(\si) \si(\ka_{1-z}(\tau)) \;.
$$

\segment{110722c}{}%PPPPPPPPPPP
Finally, we remark that the local unipotent $p$-adic Hodge morphism is an isomorphism. This follows, at least at the level of $\Qp$-points, from the fact that the morphism of mixed Tate categories from which it originates is an equivalence of categories; see segment \ref{111025b} below.

%%%%%%%%%%%%%%%%%%%%%%%%%%%%%%%%%%%%%%%
\section{The mixed Tate categories to be considered}%%%%%%%%%%%%%%%
%%%%%%%%%%%%%%%%%%%%%%%%%%%%%%%%%%%%%%%%

\segment{111025b}{}%%%%%%%%%%%%
Let $G_{\Qp}$ denote the Galois group of an algebraic closure $\QQ_p^a$ of $\Qp$ over $\Qp$.  According to Proposition 4.2 of Colmez-Fontaine \cite{ColmezFontaine},
%\mar{better ref?}
 the category $\Rep_\m{st}(G_\Qp)$ of semistable representations, is a Tannakian subcategory of the category of all representations. An object $E$ of $\Rep_\m{st}$ is crystalline if and only if the associated filtered $(\phi, N)$-module $D_\m{st}(E)$ satisfies $N=0$, and it is clear that the latter is preserved under the usual additive and Tannakian operations. Thus, the subcategory $\Rep_\m{cris}$ of crystalline representations is Tannakian. We let $\mtEtQpf(\Qp) = \langle \Qp(1)_\et \rangle_\m{mt}$ denote the mixed Tate subcategory generated by the Tate object (\ref{111026}).
 
We may consider simultaneously the category of admissible filtered $\varphi$-modules over $\Qp$. The functor $D_\m{cris}: \Rep_\m{cris}(G_\Qp) \to \tensor*[^{\m{adm}}]{\m{dR}}{}(\Qp)$ of $p$-adic Hodge theory, from crystalline representations of $G_\Qp$ to admissible filtered $\varphi$-modules, is an equivalence of tensor categories (c.f. Colmez-Fontaine \cite{ColmezFontaine}). This means, in particular, that the category of admissible filtered $\phi$-modules is Tannakian. We may thus consider the category $\mtdRQp = \langle \Qp(1)_\dR \rangle _\m{mt}$ generated by the Tate object. $D_\m{cris}$ restricts to a morphism of mixed Tate categories
\[
\mtEtQpf(\Qp) \to \mtdRQp
\,,
\]
which is an equivalence of underlying categories.

\segment{111026a}{}%%%%%%%%%%%%%%
Let $T$ denote a finite set of primes containing our special prime $p$. We consider the category of representations of $G_\QQ$ on finite dimensional $\Qp$-vector spaces, which are unramified outside of $T$, and crystalline at $p$. It is clear that the functor which takes a representation to its restriction to a given decomposition group is compatible with the additive and Tannakian operations. Moreover, being unramified is closed under the additive and Tannakian operations. So in view of segment \ref{111025b}, this category is a Tannakian subcategory of the category of all representations. We let $\mtEtQpf(\ZZ[T\inv]) = \langle \Qp(1)_\et \rangle _\m{mt}$ denote the mixed Tate subcategory generated by the Tate object. The localization morphism
\[
\mtEtQpf(\ZZ[T\inv]) \to \mtEtQpf(\Qp)
\]
determined by the choice of an embedding of algebraic closures $\QQ^a \subset \QQ_p^a$, is a morphism of mixed Tate categories.

\segment{120103e}{}%%%%%%%
We let $H_\et$ denote the Heisenberg object (\ref{111018a}) of $\mtEtQpf(\ZZ[T \inv])$ and we let $H_\dR$ denote the Heisenberg object of $\mtdRQp$. Then for $n=2$ Kim's cutter (\ref{120823-1}) may be written in our language as follows.
%ppp
\[
\xymatrix{
%1
X(\ZZ[S \inv])	\ar[r] \ar[d]_\ka				&
X(\Zp)			\ar[d]^\al					\\
\m{H}^1(H_\et)		\ar[r]_{h_2}					&
\m{H}^1(H_\dR)								}
\]
(For our purposes, it will be sufficient to regard $h_2$ as a map of pointed sets, while ignoring the underlying functors.) Indeed, this is only a matter of listing several facts established above and elsewhere. 
\begin{itemize}
\item[(i)] $\pioneet{2} = H_\et$.
\item[(ii)] $\pionedr{2} = H_\dR$.
\end{itemize}
These follow from section \ref{120129-1} above.
\begin{itemize}
\item[(iii)] Kim's $\m{H}^1_{\m f}(G_T, \pioneet{2})$ is equal to our $\m{H}^1(H_\et)$.
\end{itemize}
Indeed, both parametrize $(\Qp(0), \Qp(-1), \Qp(-2))$-framed representations of $G_T$ which are crystalline at $p$.
\begin{itemize}
\item[(iv)] For any $n$, there's a canonical bijection $\m H^1(\mtdRQp, \pionedr{n})= \pionedr{n}$; in particular, $\m H^1(H_\dR) = H_\dR$.
\end{itemize}
This follows from Lemma 3.2.3 of Hadian-Jazi \cite{Hadian}; see also segment \ref{120202-1} below for the construction of the map $\m H^1(H_\dR) \to H_\dR$.
\begin{itemize}
\item[(v)] Kim's map $\m H^1_{\m f}(G_T, \pioneet{n}) \to \pionedr{n}$ is essentially constructed at the level of $\Qp$ points, by applying the functor
\[
\mtEtQpf(\ZZ[T \inv]) \to \mtdRQp
\]
to $\pioneet{n}$-torsor-objects in $\mtEtQpf(\ZZ[T \inv])$. For $n=2$, these correspond to our Heisenberg torsor objects (\ref{111018a}).
\end{itemize}

\segment{120103f}{Claim}%%%%%%%%%%
We have $\m H^1(H_\et) = (\ZZ[S \inv]^* \otimes \Qp)^2 = \QQ_p^S \times \QQ_p^S$, and through this equality, $\ka$ sends $t \in X(\ZZ[S \inv])$ to $(t, 1-t)$. 

\medskip
\noindent
This too is merely a matter of reviewing. We recall that the apparent asymmetry is introduced when we use the embedding
\[
\xymatrix
@ R = 3mm
{
X					\ar@{^(->}[r]			&
\Gm \times \Gm						\\
t					\ar@{|->}[r]		&
(t, 1-t)								}
\]
to make the identification $\pione{1}= \QQ(1) \oplus \QQ(1)$. We then have the calculation \ref{110720d}, from which $\m H^1(\pioneet{1}) = (\ZZ[S \inv]^* \otimes \Qp)^2$. From this, by comparison with the profinite picture, it is not hard to see that the unipotent Kummer map in depth \emph{one} is given by $t \to (t, 1-t)$. As for the first equality, we have Soul\'e
vanishing
\[
\m H^1(\Qp(2)_\et) = \m H^2(\Qp(2)_\et) = 0
\]
on the one hand, and the exact sequence of segment \ref{111018f} on the other.

\segment{111026b}{Theorem}%%%%%
For each finite set $S$ of primes and $p$ a prime $\notin S$, there exist a mixed Tate category $\mtMotQ(\ZZ[S \inv])$ over $\QQ$, and a commuting triangle of mixed Tate categories
\[
\xymatrix{
\mtMotQ(\ZZ[S \inv]) \ar[d]_-{\rho^\et} \ar[dr]^-{\rho^\dR}  \\
\mtEtQpf(\ZZ[T \inv]) \ar[r] & \mtdRQp
}
\]
satisfying conditions (i)--(iv) below. Here, as usual, $T$ denotes the set $S \cup \{p\}$. 
\begin{itemize}
\item[(i)] For all $S$ and $p$, we have
\[
\m H^1(\mtMotQ(\ZZ[S \inv]), \QQ(2)) = \m H^2(\mtMotQ(\ZZ[S \inv]), \QQ(2)) = 0 \,.
\]
\item[(ii)] For all $S$ and $p$, the morphism of cohomologies
\[
\xymatrix{
\m H^1(\mtMotQ(\ZZ[S \inv]), \QQ(1)) 	\ar[r]^{\rho^\et_*} \ar@{=}[d]		&
\m H^1(\mtEtQpf(\ZZ[T \inv]), \Qp(1))	\ar@{=}[d]						\\
\ZZ[S \inv]^* \otimes \QQ					\ar@{^(->}[r]_\iota							&
\ZZ[S \inv]^* \otimes \Qp														}
\]
induced by $\rho^\et$, fits into a commuting square as shown, where $\iota$ is the canonical inclusion.
\item[(iii)] We let $H_\mot$ denote the Heisenberg object of the alleged category. For all $S$ and $p$, there's a map $\tau$ forming a commutative triangle as follows.
\[
\xymatrix{
							&
X(\ZZ[S \inv])		\ar[dl]_\tau	\ar[d]^\ka		\\
\m H^1(H_\mot)		\ar[r]_{\rho^\et_2}		&
\m H^1(H_\et)				}
\]
The composite
\[
X(\ZSinv) \xto \tau \m H^1(H_\mot) \to \m H^1(\QQ(1)^2) = (\ZSinv^* \otimes \QQ)^2
\]
sends $t \in X(\ZZ[S \inv])$ to $(t, 1-t)$.
\item[(iv)] Whenever $S'$ is a finite set of primes containing $S$ but not $p$, $\rho^\dR$ factors through a morphism of mixed Tate categories
\[
\mtMotQ(\ZZ[S\inv]) \to \mtMotQ(\ZZ[{S'}\inv])
\,.
\]
\end{itemize}

\medskip
\noindent
The proof is in paragraphs \ref{120131-4}--\ref{120131-5} below.

\segment{120131-4}{}%%%%%
We may take for $\mtMotQ(\ZZ[S\inv])$ the category of mixed Tate motives unramified outside of $S$ defined, for instance, in segment 1.7 of Deligne-Goncharov \cite{DG} (where it is denoted by $\m{MT}(\Oo_S)$). According to Proposition 1.8 of loc. cit. the $p$-adic \'etale realization $\rho^\et(E)$ of an object $E$ of $\Minv$ is unramified outside of $S$, and according to Theorem 2.2.3 of Chatzistamatiou-\"Unver \cite{ChatUnv}, $\rho^\et(E)$ is crystalline at $p$. This gives us the morphism of mixed Tate categories $\rho^\et$. We may define $\rho^\dR$ by composing with localization at $p$ (\ref{111026a}) followed by the functor (\ref{111025b}) of $p$-adic Hodge theory.

\segment{120131-5}{}
We must verify properties (i)--(iv).

(i) According to Proposition 1.9(i) of Deligne-Goncharov \cite{DelGon}, we 
\[
\m H^1(\Minv, \QQ(2)) 	= \m H^1(\Mot(\QQ), \QQ(2)) 
\,.
\]
By 1.6.11 of Deligne-Goncharov, the latter equals $ K_3(\QQ) \otimes \QQ $. Recall Borel's determination of the ranks of K-groups of a number field $F$ \cite{Borel53, Borel77}:
\[
\dim K_{2m-1}(F) \otimes \QQ = \left\{
\begin{matrix}
r+s 	& m \mbox{ odd} \\
s			& m \mbox{ even,}
\end{matrix}
\right.
\]
where $r$ denotes the number of real places and $s$ denotes half the number of complex places. In particular, for $F=\QQ$ the dimension is zero, whence the stated vanishing.

As for $\m H^2$, we have, according to Proposition 1.9(ii) of DG \cite{DelGon}, $\m H^2(E)=0$ for all $E \in \Minv$.

(ii) We have $\pioneinfty(\Gm, 1) = \pione{1}(\Gm,1) = \QQ(1)$, and $K_1(\QQ) = \QQ^*$. These isomorphisms, together with the unipotent Kummer map for $\Gm$, give rise to a commuting square as follows.
\[
\xymatrix{
\Gm(\QQ)							\ar[d] \ar@{=}[r]			&
\QQ^*								\ar[d]						\\
\m H^1(\Mot(\QQ), \QQ(1)	)		\ar@{=}[r]					&
\QQ^* \otimes \QQ											}
\]
It follows from the definition of $\Minv$ as a full subcategory of $\Mot(\QQ)$ that this square restricts to the following one.
\[
\xymatrix{
\Gm\pZSinv								\ar[d] \ar@{=}[r]			&
\ZSinv^*								\ar[d]						\\
\m H^1(\Minv, \QQ(1)	)		\ar@{=}[r]					&
\ZSinv^* \otimes \QQ											}
\]

(iii)
%
%\segment{111227e}{}%%%%%
Motivic path torsors are constructed in
%\mar{What type of a thing is a torsor in Tannakian terms? How do we know that the thing constructed by DG is such a thing?}
segment 3.12 of Deligne-Goncharov \cite{DelGon}. These give rise to the map $\tau$ in the diagram below.
\[
\xymatrix{
	% line 1
																						& 
X(\ZZ[S\inv]) 						\ar@/_/[dl]_\tau \ar@/^/[ddr]^-{\ka'} \ar@{.>}[d]^-{\ka}	\\
	% line 2
H^1(H_\mot) 						\ar@{.>}[r]^-{\rho^\et_2} \ar@/_/[drr]_-{\textit{DG}} 	&
H^1(\mtEtQpf(\ZZ[T \inv]), H_\et) 	\ar@{_{(}->}[dr]												\\
	% line 3
																						&
  																						& 
H^1(\tensor*[^{\m{mt}}]{\textrm{\'Et}}{^\Qp}(\ZZ[T\inv]), H_\et)								}
\]
We let $\tensor*[^{\m{mt}}]{\textrm{\'Et}}{^\Qp}(\ZZ[T\inv])$ denote the category of mixed Tate representations of $G_T$ which are unramified outside of $T$ (but not necessarily crystalline at $p$), and we let $\ka'$ denote the map obtained from the Kummer map, as shown. According to segment 3.13 of Deligne-Goncharov,
%\mar{This reference is weak}
 the motivic path torsors have $p$-adic \'etale realizations for $p \notin S$ which agree with the \'etale path torsors which define the map $\ka$. This gives us the map \textit{DG}, as well as the commutativity of the triangle with sides $\textit{DG}, \tau, \ka'$. According to Theorem 2.2.3 of Chatzistamatiou-Unver \cite{ChatUnv}, $p$-adic \'etale realizations of mixed Tate motives unramified at $p$ are crystalline. Thus the factorization of \textit{DG} through $\rho^\et_2$ as shown.
 
To make sense of the map $t \mapsto (t, 1-t)$, recall that the embedding of $\Gm$ in $\Gm \times \Gm$ given by the same formula, induces an isomorphism like so.
\[
\pione{1}(X, \vec{01}) \iso \pione{1}(\Gm \times \Gm, (1,1)) = \pione{1}(\Gm, 1)^2 = \QQ(1)^2
\]
Moreover, it follows from the construction of the motivic torsor map $\tau$ that the resulting square
\[
\xymatrix{
X(\QQ) 	\ar[r]					
								\ar@{^{(}->}[d] 			
															&
\m H^1(\Minv, \pione{1}) 		
								\ar@{=}[d] 					
															\\
(\Gm \times \Gm)(\QQ) 		
								\ar[r] 						
															&
(\QQ^* \otimes \QQ)^2		
															}
\]
 commutes.
 
 (iv) is a straightforward consequence of the construction of $\rho^\dR$.

%%%%%%%%%%%%%%%%%%%%%%%%%%%%%%%%%%%%%%%
\section{The de Rham splitting}%%%%%%%%%%
%%%%%%%%%%%%%%%%%%%%%%%%%%%%%%%%%%%%%%%

\segment{120202-1}{}
We denote the Heisenberg group object (\ref{111018a}) of the mixed Tate category $\DR$ by $H_\dR$. Its Lie algebra $\fk h_\dR$ has a canonical basis (given by $(dz/z)^\lor$, by $(dz/(1-z))^\lor$, and by their bracket), so the short exact sequence
\[
0 \to \Qp(2) \to \fk h_\dR \to \Qp(1)^2 \to 0
\]
of $\Qp$-vector spaces is canonically split. (This also follows from \S\ref{120129-1}.) We claim that the resulting map $\Sigma$ fits into a commuting square like so.
\[
\xymatrix{
A_2^\m{dR}					\ar[r]^-{\widetilde \Psi}					&
\Qp																		\\
\m H^1(H_\m{dR})				\ar[u]^{\la_2^\dR} \ar[r]^-\sim _-\Psi		&
H_\m{dR}						\ar[u]_\Sigma					}
\]
In it, $\la_2^\dR$ is the map $\la_2$ of section \ref{120201-1}, applied to the mixed Tate category $\DR$. To define $\Psi$, let $(E, \varphi, F)$ be a $(\Qp(0), \Qp(-1), \Qp(-2))$-framed object representing an element $[E, \varphi, F]$ of $\m H^1(H_\dR)$. Then the slope decomposition is a canonical decomposition of the underlying $\varphi$-module as
\[
E = \Qp(0) \oplus \Qp(-1) \oplus \Qp(-2)
\,.
\]
Regard $H_\dR$ as the unipotent upper triangular subgroup of $\GL(E)$:
\[
H_\dR =
\begin{pmatrix}
1 & \Qp(1)_\dR & \Qp(2)_\dR \\
& 1 & \Qp(1)_\dR \\
&& 1
\end{pmatrix}
\]
Let $F_1$ denote the direct-sum filtration on $E$, and for $u\in H_\dR$, let $F_u$ denote the filtration
\[
F^i_uE = u(F^i_1E)
\,.
\]
Then the original filtration $F$ is of the form $F= F_u$ for a unique $u \in H_\dR$, and $\Psi$ is defined by $\Psi[E,\phi,F] = u$. Thus defined, $\Psi$ is an isomorphism. This is a special case of Lemma 3.2.3 of Hadian-Jazi \cite{Hadian}.

\segment{120202-2}{}%%%%%
Now let $(E, \varphi, F)$ be a $(\Qp(0), \ast, \Qp(-2))$-framed object representing an element of $A_2^\dR$. Then $(E, \varphi)$ has slope decomposition
\[
E = \Qp(0) \oplus E_{-1} \oplus \Qp(-2)
\,,
\]
where $E_{-1}$ is isomorphic to a sum of copies of $\Qp(-1)$. Then the original filtration $F$ is uniquely of the form $F_u$ as above, except that this time $u$ is a $(\Qp$-)point of (the unipotent $\Qp$-group underlying) the unipotent group object
\[
\begin{pmatrix}
1 & E_{-1}^\lor & \Qp(2)_\dR \\
& 1 & E_{-1}(2) \\
&& 1
\end{pmatrix}
\,.
\]
Let $c$ be its component in $\Qp(2)$. Then $[E,\varphi,F] \mapsto c$ defines the map $\widetilde \Psi$. It is clear that the square commutes.

\segment{120202-3}{}%%%%%
The map $\widetilde \Psi$ may be regarded as a translation to $A_2^\dR$ of the splitting $\Sigma$ as follows.

\subsection*{Four propositions}
\begin{itemize}
\item[(1)] $A_1^\dR = \Qp$.
\item[(2)] $\Ext^2_\DR(\Qp(0), \Qp(2)) = 0$.
\item[(3)] $\Ext^1_\DR(\Qp(0), \Qp(2))= \Qp$.
\item[(4)] The map $\Ext^1_\DR(\Qp(0), \Qp(2)) \to A_2^\dR$ of segment \ref{111018f} applied to $T = \DR$ is injective and canonically split.
\end{itemize}

\begin{proof}
$A_2=\Ext^1(k(0), k(1))$ (\ref{111018d}--\ref{111018e}). So (1) and (3) follow from Proposition 1.2.7 of Chatzistamatiou-Unver \cite{ChatUnv}. (2) follows from the fact that the mixed Tate dual $U_\dR$ of $\DR$ is free; see Corollary 1.4.6 of loc. cit. Finally, to verify (4), we should elaborate on the isomorphism of (3). It is given by sending a $(\Qp(0), \Qp(2))$-framed object $(E, \varphi, F)$ to $c$, where
\[
u=
\begin{pmatrix}
1&c\\&1
\end{pmatrix}
\in
\begin{pmatrix}
1 & \Qp(2) \\ & 1
\end{pmatrix}
\]
is the element such that $F= F_u$ with respect to the slope decomposition; that is, it is given precisely by the composite
\[
\Ext^1_\DR(\Qp(0), \Qp(2)) \to A_2^\dR \xto{\widetilde \Psi} \Qp
\,,
\]
with $\widetilde \Psi$. This establishes (4).
\end{proof}

\subsection*{Corollary}
$A_2^\dR = \QQ_p^2$.

\begin{proof}
We apply the four propositions to the sequence of morphisms of segment \ref{111018f} to obtain a split exact sequence
\[
\xymatrix{
0 							\ar[r]					&
\Qp							\ar[r]					&
A_2^\dR					\ar[r] \ar@/_/[l]_{\widetilde \Psi}	&
\Qp \otimes_\Qp \Qp		\ar[r] 					&
0													}
\]
of vector spaces.
\end{proof}

%%%%%%%%%%%%%%%%%%%%%%%%%%%%%%%
\section{Linearization, evaluation at geometric points, and the twisted antisymmetry relation}%%%%%%%%%%%%%%%%%%%%
%%%%%%%%%%%%%%%%%%%%%%%%%%%%%%%

\segment{111212c}{}%%%%%%%%%%%%%%%%%%%%%
Let $S$ be a finite set of primes, and $p$ a prime $\notin S$. Then we have the mixed Tate categories $\mtdRQp$ over $\Qp$ (\ref{111025b}), $\mtMotQ(\ZZ[S \inv])$ over $\QQ$ (\ref{111026b}), and $\mtEtQpf(\ZZ[T\inv])$ over $\Qp$ (\ref{111026a}). These possess Tate objects $?(1)$ and Heisenberg objects $H_?$ (\ref{111018a}), and those possess pointed cohomology sets $\m H^1(?(1))$ and $\m H^1(H_?)$ (\ref{111018f}). Moreover, we have associated graded Hopf algebras $A^?$ (\ref{111010b}) and linearization maps $\la^?$ (\ref{111018f}). The commutative triangle of Theorem \ref{111026b} then gives rise to a commutative diagram --- the part of the figure framed by the morphisms $\rho_1^\et$, $\m{ab}^\mot$, $\la_2^\mot$, $\tensor*[^A]{\rho}{^\dR_2}$, $\la^\dR$, $\m{ab}^\dR$, $h_1$.
\begin{figure}[b]
\small	{
\[
\xymatrix
@ C=5mm @ R=10mm
{
%make source very small for tabbing to display correctly
%1
										&														&						&											& 0 \ar[d]					\\
%2
A_2^\mot \ar@{=}[dd]_-\nu \ar[rr]^-{\tensor*[^A]{\rho}{_2 ^\dR}} 	&														& A_2^\dR \ar[rr]^-{\widetilde \Psi} 	&											& \Qp(2)_\dR \ar@{.>}[d] 			\\
%3
										& \m H^1(H_\mot) \ar[dr] \ar[ul]_-{\la^\m{mot}_2} \ar@{=}[dd]^-{\m{ab}^\m{mot}} \ar[rr] 			&						& \m H^1(H_\dR) \ar[dd]^-{\m{ab}^\dR} \ar[ul]_-{\la^\dR} \ar@{=}[r]^-{\Psi}	& H_\dR \ar[dd] \ar @/_/ [u]_-{\Sigma} 	\\
%4
(A_1^\mot) ^{\otimes 2}				 			&														& \m H^1(H_\et) \ar[ur]_-{h_2} \ar@{=}[dd] 	&											&						\\
%5
\QQ^S \otimes \QQ^S \ar@{=}[u] 						& \m H^1(\QQ(1)_\mot^{\oplus2}) \ar'[r][rr] \ar[dr]_-{\rho^\et_1} \ar[ul]_-{\la^\m{mot}_1}	 	&						& \m H^1(\Qp(1)^{\oplus2}_\dR) \ar@{=}[r] 						& \Qp(1)_\dR^{\oplus 2}  \ar[d]  		\\
%6
										& \QQ^S \times \QQ^S \ar[ul] \ar[dr] \ar@{=}[u] 				& \m H^1(\Qp(1)_\et^{\oplus2}) \ar[ur]_-{h_1}  	&						& 0 		 									&						\\
%7
										&														& \QQ_p^S \times \QQ_p^S \ar@{=}[u] 	&											&						\\
%8
										& \QQ \ar@{^{(}->}[r] & \Qp											&						&											&						}
\]
		}
\end{figure}
If we apply paragraph \ref{111018f} to the mixed Tate category $\mtMotQ(\ZZ[S \inv])$, we obtain the square formed by the morphisms $\nu$, $\la_2^\mot$, $\m{ab}^\mot$, $\la^\mot_1$. The vanishing results \ref{111026b}(i) imply that $\nu$ and $\m{ab}^\mot$ are iso, as shown. Segments \ref{120202-1}--\ref{120202-2} give us the map $\widetilde{\Psi}$, as well as the commutativity of the square formed by $\widetilde{\Psi}$, $\Sigma$, $\Psi$, $\la^\dR$. Also included in the diagram, is the short exact sequence defining $H_\dR$ as a central extension on the right. (One of its arrows --- the dotted one, doesn't commute with the rest of the diagram.) To complete the construction of the diagram, we recall that
\[
\m H^1(\mtEtQpf(\ZZ[T \inv], \Qp(1)) = \m H^1_\m f(G_T, \Qp(1))
\]
(\ref{120103e}), so the equality near the bottom follows from Corollary \ref{110720d}, and the two equalities that run parallel to this one follow from Theorem \ref{111026b} (ii), and from paragraphs \ref{111018d}--\ref{111018e}, respectively.

\subsection*
%{111026e}
{Theorem}%TTTTTTTTT
The third coordinate $h_{1,3}$ of the global unipotent $p$-adic Hodge morphism for the thrice punctured line in depth two over $\ZZ[S \inv]$ is bilinear.

\begin{proof}%pppppppppp
We have $h_{1,3} = \Sigma \circ \Psi \circ h_2 $. Since these maps come from morphisms of $\Qp$-schemes by taking $\Qp$-points, we may check bilinearity after restricting to the $\QQ$-lattice $\m H^1(\QQ(1)^{\oplus2}_\mot)$. By the commutativity of the diagram, we may thus replace the composite $\Sigma \circ \Psi \circ h_2$ by $\widetilde{\Psi} \circ \tensor*[^A]{\rho}{^\dR_2} \circ \nu \inv \circ \la_1^\mot$. Of these, $\la_1^\mot$ is $\QQ$-bilinear, $\nu$ is $\QQ$-linear, $\tensor*[^A]{\rho}{^\dR_2}$ is linear over the embedding $\QQ \subset \Qp$, and $\widetilde \Psi$ is $\Qp$-linear.
\end{proof}

\segment{111227a}{}%%%%%%%%%%%%%%%%%%
We let $\tensor*[_\QQ]{h}{_{1,3}}$ denote the linear map $\QQ^S \otimes \QQ^S \to \Qp(2)_\dR$ (as well as the associated bilinear map) defined by the diagram of segment \ref{111212c}. We remark that it follows from Theorem \ref{111026b}(iii) that formation of $\tensor*[_\QQ]{h}{_{1,3}}$ is compatible with change of $S$: if $S'$ is a finite set of primes containing $S$ but not $p$, then $\tensor*[_\QQ]{h}{_{1,3}}(S)$ factors through $\tensor*[_\QQ]{h}{_{1,3}}(S')$. If $t \in \Gm\pZSinv$, we write $t$ again for its image in $\QQ^S$. 

\segment{111227c}{Proposition}%PPPPPggg
For any $t \in X(\ZZ[S \inv])$, we have
\[_\QQ h_{1,3}(t, 1-t) = - \Li (t)
\,.
\]

%\medskip
%\noindent
%The proof follows in segments \ref{111227e}--\ref{111227f} below.

%\ssegment{111227f}{}%%%%%%%%%%%%%%
\begin{proof}
We combine Kim's cutter (\ref{120103e}) and the triangle of segment \ref{111026b}(iii) with the isomorphisms $\m H^1(H_\mot) \xto{\sim} \m H^1(\QQ(1)^2)$ (\ref{111212c}) and $\m H^1(\QQ(1)) $  $ = \QQ^S$ (\ref{111026b}(ii)), to obtain the commutative diagram below.
\[
\xymatrix{
	% Line 1
																			&
X(\ZZ[S \inv])			\ar[dl]_-\tau \ar[d]^-\ka \ar[r]							&
X(\Zp)					\ar[d]												\\
	% Line 2
\m H^1(H_\mot)			\ar@{=}[d] \ar[r]_{\rho^\et_2}							&
\m H^1(H_\et)				\ar[r]												&
\m H^1(H_\dR)			\ar[r]												&
\Qp(2)_\dR																	\\
	% Line 3
\QQ^S \times \QQ^S														}
\]
In this diagram, a point $t \in X(\ZZ[S \inv])$ corresponds to $(t, 1-t) \in \QQ^S \times \QQ^S$ (Theorem \ref{111026b}(iii)), and to $-\Li(t) \in \Qp(2)_\dR$ (Proposition \ref{110721c}). This completes the proof.
\end{proof}

\hyphenation{an-ti-sym-me-try}

\segment{111227d}{Proposition}%PPPPPttt
For any $u, v \in \Gm(\ZZ[S \inv])$, we have the twisted-anti\-symmetry relation
\[
_\QQ h_{1,3}(u, v) + \tensor*[_\QQ] {h}{_{1,3}} (v,u) = (\log u)(\log v)
\,.
\]

\begin{proof}%pppppp
If we continue to write ``$\cdot$" for the composition law of the $\QQ$-vector space $A^\mot_1 \otimes A^\mot_1 = \QQ^S \otimes \QQ^S$, then we must choose a different notation for the multiplication laws of $A^\mot$, $A^\dR$: let it be ``$*$". Let $^A \rho ^\dR$ denote the map $A^\mot \to A^\dR$ induced by $\rho^\dR$. By Corollary \ref{110720d}, we have $^A\rho^\dR (u) = \log u$, and similarly for $v$. By Proposition \ref{120103c},
\[
\nu \inv \big((u \otimes v)\cdot(v \otimes u)\big) = u * v
\,.
\]
Combining these observations, we obtain
\[
^A\rho^\dR \Big(  \nu \inv \big( (u \otimes v)\cdot(v \otimes u) \big) \Big) = ( \log u ) * ( \log v )
\,,
\]
so we will be done as soon as we've shown that for $z, z' \in A_1^\dR = \Qp$, we have $\widetilde \Psi(z*z') = zz' $ (multiplication in $\Qp$).

A $(\Qp(0), \Qp(-1))$-framed object $(E, \varphi, F)$ which corresponds to $z$ is given by $(E, \varphi) = \Qp(0) \oplus \Qp (-1)$, and (in the notation of segment \ref{120202-1}) by $F = F_u$, where
\[
u =
\begin{pmatrix}
1 	& z \\
	& 1
\end{pmatrix}
\in
\begin{pmatrix}
1	& \Qp(1)	\\
	& 1
\end{pmatrix}
\,.
\]
We write $E'$, etc. when we replace $z$ by $z'$. According to Proposition \ref{120103b}, the product, in terms of framed objects, is given by $E \varoast E'$. This, in turn, is the framed object $(E'', \varphi'', F'')$ given by
\[
(E'', \varphi'') = (E, \varphi) \otimes (E', \varphi') = \Qp(0) \oplus \Qp(-1)^2 \oplus \Qp(-2)
\,,
\]
and by $F'' = F_{u''}$, where $u'' = u \otimes u'$. Recalling (from segment \ref{120202-2}) the construction of $\widetilde \Psi$, we find that $\widetilde \Psi(E'') = zz'$, indeed.
\end{proof}

%%%%%%%%%%%%%%%%%%%%%%%%%%%%%%%%%%%%%%%
\section{Vanishing of $K_2 \otimes \QQ$ and determination of coefficients}%%
\label{Determination}%%%%%%%%%%%%%%%%%%%%%%%%%%%%%
%%%%%%%%%%%%%%%%%%%%%%%%%%%%%%%%%%%%%%%

\segment{111211a}{}%%%%%%%%%%%%%%%%%%%%%%%%
We begin with a construction closely related to Tate's computation of the Milnor $K$-group $K_2(\QQ)$ (c.f. \S11 of Milnor \cite{milnor}). Let $E = \QQ \otimes \QQ^*$, which we regard as a $\QQ$-vector space, written multiplicatively. So a linear combination looks like so.
\[
v_1^{a_1} \cdots v_n^{a_n}
\]
We commit a standard abuse of notation, denoting the element $1 \otimes a$, for $a \in \QQ^*$, simply by $a$. (In this notation, we have $-1 = 1$.) We define two subsets $R_s$ and $R_d$ of $E \otimes E$:
\[
R_s = \big\{ (a \otimes b) \cdot (b \otimes a) \; \big| \; a, b \in \QQ^* \big\}
\]
\[
R_d = \big\{ a \otimes b \; \big| \; a, b \in \QQ^*, \; a+b=1 \big\} \;,
\]
and we let $R = R_s \cup R_d$. Our construction gives an algorithm for expressing an arbitrary generator $q \otimes q'$ of $E \otimes E$ as a linear combination
\[
q \otimes q' = \prod_k \big( (u_k \otimes v_k) \cdot (v_k \otimes u_k) \big)^{s_k} \cdot \prod_l \big( t_l\otimes (1-t_l) \big)^{d_l}
\tag{$*$}
\]
of elements of $R$, with
\[
|u_k|, |v_k|, |t_l|, |1-t_l| \le \max\{|q|, |q'| \}
\]
for all $k$ and $l$.

\ssegment{0}{}%%%%%%%%%%%%%%%%%%
It is clear how to write an arbitrary element as a linear combination of elements of $R_s$ and elements of the form $q \otimes p$, with $0 < q < p $ both prime, so we may restrict our attention to a generator of this form. We show how to write such a generator as a linear combination of elements of $R_d$ and elements $a \otimes b$ with $|a|, |b| < p$. Iteration of this procedure then produces the desired expression.

\ssegment{111211b}{}%%%%%%%%%%%%%%%%%%%%%%%%
Define integers $0 < z_i < p$, $0 \le r_i < p$ recursively by
\begin{align*}
z_0 &= q  \\
qz_i &= z_{i+1} + p r_{i+1} \tag{$*$}
\end{align*}
and define $f_i \in E \otimes E$, for $i \in \NN$, by
\[
f_i = \left\{
\begin{matrix}
(\frac{z_i}{qz_{i-1}} \otimes p)\inv & \mbox{if } r_i \neq 0 \\
1 & \mbox{if } r_i=0 \;.
\end{matrix}
\right.
\]

\subsection*{Lemma}%LLLLLLLLLLLLLL
For each integer $m \ge 0$,
$
q^{m+1} \otimes p = (z_m \otimes p)f_m \cdots f_1
\,.
$

\begin{proof}%ppppppppppppppp
The base case reads $q \otimes p = z_0 \otimes p$, which is true. Now consider an arbitrary $m$ for which the lemma holds. There are two cases to consider, depending on whether $r_{m+1}$ does or does not $=0$. Either way, we have
\begin{align*}
(z_{m+1} \otimes p)f_{m+1} \cdots f_1 
&= (qz_m \otimes p)f_m \cdots f_1 \\
&= (q \otimes p)(z_m \otimes p)f_m \cdots f_1 \\
&= (q\otimes p)(q^{m+1} \otimes p) \\
&= q^{m+2} \otimes p
\,. \qedhere
\end{align*}
\end{proof}

\ssegment{111211c}{}%%%%%%%%%%%%%%%%%%%%%%%
Reducing modulo p, equation \ref{111211b}($*$) implies $z_i \equiv q^{i+1}$ for each $i$. So if $\al$ is the order of $q$ in $\FF_p^*$, we have $z_{\al-1}=1$. Hence, by Lemma \ref{111211b},
\[
q \otimes p = (f_{\al-1} \cdots f_1)^{1/\al}
\,.
\]
We've written our generator as a linear combination of $f_i$'s. Moreover, for each $i$,
\[
f_i = \left( \frac{z_i}{qz_{i-1}}  \otimes \frac{r_i}{qz_{i-1}} \right) \left( \frac{z_i}{qz_{i-1}} \otimes \frac{pr_i}{qz_{i-1}} \right)\inv
\]
is a product of two terms. The first is a tensor product of nonzero rational numbers of absolute value less than $p$; the second is (a power of) an element of $R_d$. This completes the construction.

\segment{111212a}{}%%%%%%%%%%%%%%%%%%%%%%%%%%
For each ordered pair $(q,q')$ of prime numbers, and $p$ a prime $> q, q'$, we let $\tensor*[_\QQ ]{ h}{_{1,3}^{q,q'}}$ denote the $p$-adic number given in the notation of segment \ref{111211a} by
\[
\tensor*[_\QQ ] {h}{^{q,q'}_{1,3}} = \sum_k s_k (\log u_k)(\log v_k) -\sum_l d_l \Li(t_l) \,.
\]
Fix a finite set $S$ of primes, and $p$ a prime greater than the primes of $S$. We recall that $h_{1,3}$ denotes the third coordinate of the unipotent $p$-adic Hodge morphism associated to $S$, and that $_\QQ h_{1,3}$ is its restriction to $\QQ^S \times \QQ^S$. Let $q_1, \dots, q_s$ denote the elements of $S$. We use the notation $q^x=q_1^{x_1} \cdots q_s^{x_s}$ for an element of $\QQ^S$.

%\ssegment{111227g}
\subsection*{Theorem}%TTTTTTTTTTT
For any pair $q^x$, $q^y$, of elements of $\QQ^S$, we have the formula
\[
_\QQ h_{1,3}(q^x, q^y)
=
\sum_{1 \le i,j \le s} \tensor*[_\QQ ] {h}{^{q_i, q_j}_{1,3}}x_i y_j
\,.
\]

\begin{proof}%ppppppppp
We've seen that $_\QQ h_{1,3}$ is bilinear (\ref{111212c}). So after fixing integers $i, j \in [1,s]$, it remains to show that
\[
_\QQ h_{1,3}(q_i, q_j) = \tensor*[_\QQ ] {h}{^{q_i, q_j}_{1,3}} \,.
\]
By compatibility of formation of $_\QQ h_{1,3}$ with change of $S$ (\ref{111227a}), after possibly replacing $S$ by a larger set of primes $<p$, we may assume that equation \ref{111211a}($*$) holds in $\QQ^S \otimes \QQ^S$. The theorem thus follows from Propositions \ref{111227c} and \ref{111227d} and Theorem \ref{111212c}.
\end{proof}

\segment{120202-4}{}%%%%%%%
We retain the notation of segment \ref{111212a}.

\subsection*{Theorem}
The map $\tensor*[^A]{\rho}{^\dR_2}:A_2^\mot \to A_2^\dR$ induced by the de Rham realization functor $\rho^\dR: \Minv \to \DR$ is given, through the isomorphisms $A_2^\mot = \QQ^S \otimes \QQ^S$ and $A_2^\dR = \Qp \times \Qp$ (\ref{120202-3}), by the formula
\[
q^x \otimes q^y \mapsto \big( \sum_{i,j=1}^s (\log q_i)(\log q_j) x_iy_j , \sum_{i,j=1}^s \tensor*[_\QQ ] {h}{_{1,3}^{q_i, q_j}} x_i y_j  \big)
\,.
\]

\begin{proof}
This is immediate from Theorem \ref{111212a} in view of the construction of the isomorphism $A_2^\dR = \QQ_p^2$.
\end{proof}

%%% Local Variables: 
%%% mode: latex
%%% TeX-master: "CKTwo12"
%%% End: 
\section{Examples} \label{Examples}

%%%%%%%%%%%%%%%%%%%%%
\segment{Examples1}{Siegel's theorem}

There are essentially two case where our results produce an {\em explicit}
Coleman function vanishing on all integral points of
$X=\PP^1\setminus\{0,1,\infty\}$. Suppose first that $S=\{\ell\}$ consists of a single
prime $\ell$. Then 
\[
   {\rm rank}\;\ZZ[S^{-1}]^* = 1
\]
and therefore
\[
    \dim_{\QQ_p} \m H^1_\m{f}(G_T,\pioneet{2}) = 2 < 
       \dim_{\QQ_p} \pionedr{2} = 3.
\]
It follows that the image of the map $h_2$ is not Zariski dense, and we can
conclude that the set $X(\ZZ[S^{-1}])$ is finite.\footnote{Actually, for
  $\ell\neq 2$ the set $X(\ZZ[S^{-1}])$ is empty for trivial reasons.}

For the second case, let $K$ be a real quadratic number field, and let
$\mathcal{O}_K$ denote the ring of integers of $K$. We set $S:=\emptyset$ and
$T:=\{p\}$. Then we also have
\[
   \dim_{\QQ_p} \m H^1_\m{f}(G_T,\QQ_p(1))={\rm rank}\,\mathcal{O}_K^\times=1
\]
and 
\[
   \dim_{\QQ_p} \m H^1_\m{f}(G_T,\QQ_p(2))=0,
\]
by Soul\'e's vanishing result. If we choose the prime $p$ such that it splits in
$K$ then the map
\[
      h_2: \m H^1_\m{f}(G_T,\pioneet{2})\to \pionedr{2}
\]
is defined in the same way as before, and the dimensions of source and target
are two and three, respectively. It follows as before that $X(\mathcal{O}_K)$
is finite.

Somewhat surprisingly, there exists a {\em single} Coleman function vanishing
simultaneously on all integral points of all these examples. This follows from
our results as follows. Let $R$ denote the ring $\ZZ[\ell^{-1}]$ or
$\mathcal{O}_K$ and let $\epsilon$ be a generator of the free part of
$R^\times$. For instance, if $R=\ZZ[\ell^{-1}]$ then we can set
$\epsilon:=\ell$. Note that the set $X(R)$ can be identified with solutions of
the {\em unit equation} 
\[
       \pm \epsilon^x \pm \epsilon^y = 1, \qquad x,y\in\ZZ.
\]

We identify $\m H^1_\m{f}(G_T,\pioneet{2})$ with affine two-space
$\mathbb{A}_{\QQ_p}^2$, by writing an arbitrary element in the form
$(\epsilon^x,\epsilon^y)$, with $x,y\in\QQ_p$. The map
\[
    \kappa: X(R) \to \m H^1_\m{f}(G_T,\pioneet{2})
\]
is then given by the rule
\[
      \kappa(z) = (\epsilon^x,\epsilon^y), \quad
        \text{with $z=\pm \epsilon^x$, $1-z=\pm \epsilon^y$.}
\]
Also, we identify $\pionedr{2}$ with affine three-space in such a way that the
map $\alpha:X(\ZZ_p)\to\pionedr{2}$ is given by
\[
     \alpha(z)=(\log(z),\log(1-z),-\Li(z))).
\]
Our main result says that the third coordinate of $h_2$ is bilinear, and hence
\[
      h_2(\epsilon^x,\epsilon^y) = (x\log(\epsilon),y\log(\epsilon),cxy),
\]
for some constant $c\in\QQ_p$. The {\em twisted anti-symmetry relation}
(Proposition \ref{111227d}) shows that $c=\log(\epsilon)^2/2$. We conclude that the
image of $h_2$ is contained (and in fact equal to) the subvariety given by the
equation 
\[
     2w= uv.
\]
Therefore, by the commutativity of the Kim's cutter (\ref{120823-1}), the Coleman
function
\[
     D(z) := 2\Li(z)-\log(z)\log(1-z)
\]
vanishes on the subset $X(R)\subset X(\ZZ_p)$. As announced, this function is
independent of $R$ (but does depend, of course, on $p$). 

By explicit computations, using Sage and \cite{Lip}, we found that for $p=11$
the function $D(z)$ has exactly $9$ distinct zeroes on $X(\ZZ_{11})$. Moreover, by
approximating theses zeroes with sufficient precision, it was easy to `guess'
that these are 
\[
   z=-1,2,1/2,\frac{1\pm\sqrt{5}}{2},\frac{-1\pm\sqrt{5}}{2},
       \frac{3\pm\sqrt{5}}{2}.
\]
Now the first three of these points are $\ZZ[1/2]$-integral points on $X$, and
the latter six are $\mathcal{O}_K$-integral points, for
$K=\QQ(\sqrt{5})$. Hence our guess must be correct, and we have proved the
following (well known) result.

\subsection*{Proposition} Let $R$ be either $\ZZ[\ell^{-1}]$, for a prime
$\ell$, or $R=\mathcal{O}_K$ for a real quadratic number field $K$. Suppose,
moreover, that $p=11$ splits in $K$. Then the set $X(R)$ is empty, except if
$R=\ZZ[1/2]$ or $R=\ZZ[(1+\sqrt{5})/2]$, in which case we have
  \[
       X(\ZZ[1/2]) = \{-1,2,1/2\}
  \]
  or
  \[
       X(\ZZ[(1+\sqrt{5})/2]) = \{\frac{1\pm\sqrt{5}}{2},
            \frac{-1\pm\sqrt{5}}{2},\frac{3\pm\sqrt{5}}{2}\}.
  \]

Actually, the restriction that $p=11$ splits in $K$ in unnecessary, see e.g.\
\cite{Nagell}. It is not clear how to get rid of this assumption within the
context of our results. The problem is that for primes $p>11$ the function
$D(z)$ seems to always have more than $9$ zeroes, and if one of them is a
transcendental number, then one can approximate it to arbitrary precision, but
one does not know how to prove that it does not correspond to an $R$-integral
point.  

%%%%%%%%%%%%%%%%%%%%%%%%%%%
\segment{Examples2}{Exponential Diophantine equations in two variables}

If $|S|>1$ then our results do not produce an explicit Coleman function
vanishing on $X(\ZZ[S^{-1}])$. However, we can find such functions vanishing
on certain subsets corresponding to solutions of exponential diophantine
equations of the form
\begin{equation} \label{exeq1}
    a\cdot b^x-c\cdot d^y = 1,
\end{equation}
with fixed $a,b,c,d\in\QQ^\times$. Namely, if $S$ contains all primes which
occur in the prime decomposition of $a,b,c,d$, then for any solution
$(x,y)\in\ZZ^2$ of \eqref{exeq1}, $z=ab^x$ is a $\ZZ[S^{-1}]$-integral point
of $X$. Moreover, the image of $z=ab^x$ under the map $\kappa$ lies in a
linear subspace $V\subset \m H^1_\m{f}(G_T,\pioneet{2})$ of dimension $2$, where
$V$ is determined by \eqref{exeq1}. Since $h_2(V)\subset\pionedr{2}$ is a
proper algebraic subspace for dimension reasons, we can find a Coleman
function $f(z)$ in $X(\ZZ_p)$ vanishing on all solutions of
\eqref{exeq1}. This proves finiteness of the set of solutions of
\eqref{exeq1}. However, it does not give an algorithm to determine the
solution set. 

We illustrate this by the following example. Consider the exponential
Diophantine equation 
\begin{equation} \label{exeq2}
    7^x-3\cdot 2^y=1.
\end{equation}
We set $S:=\{2,3,7\}$ and choose a prime $p\not\in S$. We identify
$\m H^1_\m{f}(G_T,\pioneet{2})$ with $\mathbb{A}^6_{\QQ_p}$, using coordinates
$x_1,x_2,x_3,y_1,y_2,y_3$ and writing an arbitrary element as
$(2^{x_1},3^{x_2},7^{x_3},2^{y_1},3^{y_2},7^{y_3})$. Then the image of a point
$z=7^x$ corresponding to a solution of \eqref{exeq2} under $\kappa$ lies in
the linear subspace
\[
    V = \{(0,0,x,y,1,0) \mid x,y\in\QQ_p\} \subset \m H^1_\m{f}(G_T,\pioneet{2}).
\]
The third coordinate $h_{1,3}$ of the map
$h_2: \m H^1_\m{f}(G_T,\pioneet{2})\to\pionedr{2}$ is bilinear; thus its restriction
to $V$ is of the form 
\[
    V\cong\mathbb{A}^2 \to \QQ_p, \quad (x,y) \mapsto c_{7,2}xy + c_{7,3}x,
\]
where $c_{7,2}=h_2(7,2)$, $c_{7,3}=h_2(7,3)$. It follows that the image
$h_2(V)\subset \pionedr{2}\cong\mathbb{A}^3$ is contained in the hypersurface
with equation 
\[
   u\cdot(c_{7,2}\cdot v-c_{7,2}\log(3)+c_{7,3}\log(2)) = \log(2)\log(7)\cdot w.
\]
We conclude, as before, that the Coleman function
\[
    f(z):=\log(2)\log(7)\cdot \Li(z) +\log(z)\cdot
       (c_{7,2}\cdot \log(1-z)-c_{7,2}\log(3)+c_{7,3}\log(2))
\]
vanishes on the points of $X(\ZZ_p)$ corresponding to solutions of
\eqref{exeq2}. 

It remains to compute the coefficients $c_{7,3},c_{7,2}$. We follow the
algorithm offered in Section \ref{Determination}. We obtain
\[
    7\otimes 2 = (7\otimes 8)^{1/3}=(-7\otimes(1+7))^{1/3}
\]
and
\[
\begin{split}
     7\otimes 3 & = (7\otimes 6)\cdot(7\otimes 2)^{-1} \\
          & = (7\otimes(1-7))\cdot(-7\otimes(1+7))^{1/3},
\end{split}
\]
and conclude that
\[
   c_{2,7}=-\frac{1}{3}\Li(-7), \quad 
   c_{3,7}=-\Li(7)+\frac{1}{3}\Li(-7).
\]

Now we choose $p=5$. It turns out that $f(z)$ has exactly $6$ zeroes on
$X(\ZZ_5)$. Of these, three can be guessed to be
\[
       z=-1,7,49.
\]
The values $z=7,49$ correspond to solutions of \eqref{exeq2}, namely
$(x,y)=(1,1)$ and $(x,y)=(2,4)$, and therefore they are indeed zeroes of
$f$. The value $z=-1$ does not correspond to a solution, but it can easily be
shown to be a zero of $f(z)$, using the identity $D(-1)=0$ established in the
previous subsection.

The other three zeroes of $f(z)$ in $X(\ZZ_p)$ do not seem to be well
approximated by algebraic numbers of small height, so our guess is that they
are transcendental. But it is not clear how to prove this with our
methods. 

It is a nontrivial exercise in elementary number theory to show that
\eqref{exeq2} has exactly the two solutions $(x,y)=(1,1),(2,4)$. So far, we
are not able to prove this using the nonabelian Chabauty method. All we get is
that there are at most five solutions. 

\bibliography{references}

\begin{thebibliography}{BDCKW}

\bibitem[BD]{beilinson-deligne}
A.~Be{\u\i}linson and P.~Deligne.
\newblock Interpr\'etation motivique de la conjecture de {Z}agier reliant
  polylogarithmes et r\'egulateurs.
\newblock In {\em Motives ({S}eattle, {WA}, 1991)}, volume~55 of {\em Proc.
  Sympos. Pure Math.}, pages 97--121. Amer. Math. Soc., Providence, RI, 1994.

\bibitem[BDCKW]{nabsd}
J.~Balakrishnan, I.~Dan-Cohen, M.~Kim, and S.~Wewers.
\newblock A non-abelian conjecture of {B}irch and {S}winnerton-{D}yer type for
  hyperbolic curves.
\newblock Preprint. arXiv:1209.0640v1.

\bibitem[BdJ]{Lip}
Amnon Besser and Rob de~Jeu.
\newblock {${\rm Li}^{(p)}$}-service? {A}n algorithm for computing {$p$}-adic
  polylogarithms.
\newblock {\em Math. Comp.}, 77(262):1105--1134, 2008.

\bibitem[Ber]{Berger}
Laurent Berger.
\newblock An introduction to the theory of {$p$}-adic representations.
\newblock In {\em Geometric aspects of {D}work theory. {V}ol. {I}, {II}}, pages
  255--292. Walter de Gruyter GmbH \& Co. KG, Berlin, 2004.

\bibitem[Bes]{Besser}
Amnon Besser.
\newblock Coleman integration using the {T}annakian formalism.
\newblock {\em Math. Ann.}, 322(1):19--48, 2002.

\bibitem[BGSV]{BGSV}
A.~A. Be{\u\i}linson, A.~B. Goncharov, V.~V. Schechtman, and A.~N. Varchenko.
\newblock Aomoto dilogarithms, mixed {H}odge structures and motivic cohomology
  of pairs of triangles on the plane.
\newblock In {\em The {G}rothendieck {F}estschrift, {V}ol.\ {I}}, volume~86 of
  {\em Progr. Math.}, pages 135--172. Birkh\"auser Boston, Boston, MA, 1990.

\bibitem[Bor1]{Borel53}
Armand Borel.
\newblock Sur la cohomologie des espaces fibr\'es principaux et des espaces
  homog\`enes de groupes de {L}ie compacts.
\newblock {\em Ann. of Math. (2)}, 57:115--207, 1953.

\bibitem[Bor2]{Borel77}
Armand Borel.
\newblock Cohomologie de {${\rm SL}_{n}$} et valeurs de fonctions zeta aux
  points entiers.
\newblock {\em Ann. Scuola Norm. Sup. Pisa Cl. Sci. (4)}, 4(4):613--636, 1977.

\bibitem[CF]{ColmezFontaine}
Pierre Colmez and Jean-Marc Fontaine.
\newblock Construction des repr\'esentations {$p$}-adiques semi-stables.
\newblock {\em Invent. Math.}, 140(1):1--43, 2000.

\bibitem[Col]{Coleman}
Robert~F. Coleman.
\newblock Dilogarithms, regulators and {$p$}-adic {$L$}-functions.
\newblock {\em Invent. Math.}, 69(2):171--208, 1982.

\bibitem[CU]{ChatUnv}
A.~Chatzistamatiou and S.~\"Unver.
\newblock On p-adic periods for mixed tate motives over a number field.
\newblock arXiv:1110.0923v1 [math.AG].

\bibitem[DCW]{HCE}
I.~Dan-Cohen and S.~Wewers.
\newblock The {H}eisenberg coboundary equation: appendix to \textit{Explicit
  Chabauty-Kim theory}.
\newblock Unpublished note, available from
  \url{https://www.uni-due.de/~hm0146/}.

\bibitem[Del]{Deligne89}
P.~Deligne.
\newblock Le groupe fondamental de la droite projective moins trois points.
\newblock In {\em Galois groups over {${\bf Q}$} ({B}erkeley, {CA}, 1987)},
  volume~16 of {\em Math. Sci. Res. Inst. Publ.}, pages 79--297. Springer, New
  York, 1989.

\bibitem[DG1]{DelGon}
Pierre Deligne and Alexander~B. Goncharov.
\newblock Groupes fondamentaux motiviques de {T}ate mixte.
\newblock {\em Ann. Sci. \'Ecole Norm. Sup. (4)}, 38(1):1--56, 2005.

\bibitem[DG2]{DG}
Michel Demazure and Pierre Gabriel.
\newblock {\em Groupes alg\'ebriques. {T}ome {I}: {G}\'eom\'etrie alg\'ebrique,
  g\'en\'eralit\'es, groupes commutatifs}.
\newblock Masson \& Cie, \'Editeur, Paris, 1970.
\newblock Avec un appendice {{\i}t Corps de classes local} par Michiel
  Hazewinkel.

\bibitem[Fur1]{FurushoI}
Hidekazu Furusho.
\newblock {$p$}-adic multiple zeta values. {I}. {$p$}-adic multiple
  polylogarithms and the {$p$}-adic {KZ} equation.
\newblock {\em Invent. Math.}, 155(2):253--286, 2004.

\bibitem[Fur2]{FurushoII}
Hidekazu Furusho.
\newblock {$p$}-adic multiple zeta values. {II}. {T}annakian interpretations.
\newblock {\em Amer. J. Math.}, 129(4):1105--1144, 2007.

\bibitem[Gon]{GonMPMTM}
A.~B. Goncharov.
\newblock Multiple polylogarithms and mixed {T}ate motives.
\newblock arXiv:math/0103059v4 [math.AG].

\bibitem[Hai]{Hain92}
Richard~M. Hain.
\newblock Classical polylogarithms.
\newblock In {\em Motives ({S}eattle, {WA}, 1991)}, volume~55 of {\em Proc.
  Sympos. Pure Math.}, pages 3--42. Amer. Math. Soc., Providence, RI, 1994.

\bibitem[HJ]{Hadian}
Majid Hadian-Jazi.
\newblock Motivic fundamental groups and integral points.
\newblock \url{http://hss.ulb.uni-bonn.de/2010/2217/2217.pdf}, 2010.

\bibitem[Kim1]{kimi}
Minhyong Kim.
\newblock The motivic fundamental group of {$\mathbb P^1\setminus
  \{0,1,\infty\}$} and the theorem of {S}iegel.
\newblock {\em Invent. Math.}, 161(3):629--656, 2005.

\bibitem[Kim2]{kimii}
Minhyong Kim.
\newblock The unipotent {A}lbanese map and {S}elmer varieties for curves.
\newblock {\em Publ. Res. Inst. Math. Sci.}, 45(1):89--133, 2009.

\bibitem[Mil]{milnor}
John Milnor.
\newblock {\em Introduction to algebraic {$K$}-theory}.
\newblock Princeton University Press, Princeton, N.J., 1971.
\newblock Annals of Mathematics Studies, No. 72.

\bibitem[Nag]{Nagell}
Trygve Nagell.
\newblock Sur une propri\'et\'e des unit\'es d'un corps alg\'ebrique.
\newblock {\em Ark. Mat.}, 5:343--356 (1964), 1964.

\bibitem[Ols]{OlssonTowards}
Martin~C. Olsson.
\newblock Towards non-abelian {$p$}-adic {H}odge theory in the good reduction
  case.
\newblock {\em Mem. Amer. Math. Soc.}, 210(990):vi+157, 2011.

\bibitem[Ser]{corloc}
Jean-Pierre Serre.
\newblock {\em Corps locaux}.
\newblock Hermann, Paris, 1968.
\newblock Deuxi{\`e}me {\'e}dition, Publications de l'Universit{\'e} de
  Nancago, No. VIII.

\bibitem[Tat]{Tate}
John Tate.
\newblock Relations between {$K_{2}$} and {G}alois cohomology.
\newblock {\em Invent. Math.}, 36:257--274, 1976.

\bibitem[Vol]{Vologodsky}
Vadim Vologodsky.
\newblock Hodge structure on the fundamental group and its application to
  {$p$}-adic integration.
\newblock {\em Mosc. Math. J.}, 3(1):205--247, 260, 2003.

\end{thebibliography}

\bibliographystyle{alphanum}

\bigskip

\Small\textsc{I.D.: Fakult\"at f\"ur Mathematik, 
Universit\"at Duisburg-Essen,
Universit\"atsstrasse 2, 
45117 Essen, 
Germany}

\Small\textsc{S.W.: Universit\"at Ulm,
Institut f\"ur Reine Mathematik,
Helmholtzstrasse 18,
89069 Ulm, Germany}

\smallskip

\textit{E-mail address (I.D.):} \texttt{ishaidc@gmail.com}

\smallskip

\textit{E-mail address (S.W.):} \texttt{stefan.wewers@uni-ulm.de}

\end{document}